\documentclass[11pt]{article}

\usepackage{amssymb,amsmath,amsthm}

\usepackage{graphicx, tikz}

\usepackage{color}
\definecolor{darkgreen}{rgb}{0,.4,0.2}
\definecolor{darkagenta}{rgb}{.5,0,.5}
\definecolor{darkred}{rgb}{0.85,0,0}
\definecolor{darkblue}{rgb}{0,0,.6}
\definecolor{lightgray}{gray}{.95}
\definecolor{rgrey}{rgb}{.8,0.4,.4}  
\definecolor{grey}{rgb}{.5,.5,.5}  

\newcommand{\grijs}[1]{{\color{grey}#1}}

\newcommand{\N}{{\mathbb{N}}}

\newcommand{\A}{{\rm A}}
\newcommand{\B}{{\rm B}}
\newcommand{\C}{{\rm C}}
\newcommand{\D}{{\rm D}}
\newcommand{\F}{\rm \scriptstyle F}
\newcommand{\G}{\rm \scriptstyle G}
\newcommand{\Ha}{\rm \scriptstyle H}

\newcommand{\duoZ}[2]{\fbox{$\begin{aligned} w & #1\\[-.1cm] R_{\cdot\!w} &  #2  \end{aligned}$}}
\newcommand{\duoO}[2]{\fbox{$\begin{aligned} w\; &\;  #1\\[-.1cm] R_w &  #2  \end{aligned}$}}
\newcommand{\mor}[2]{\Big\{\begin{aligned} a & \rightarrow #1\\[-.1cm] b & \rightarrow #2  \end{aligned}}
\newcommand{\xF}{x_{\rm \scriptstyle F}}
\newcommand{\xG}{x_{\rm \scriptstyle G}}
\newcommand{\xH}{x_{\rm \scriptstyle H}}

\newtheorem*{theorem*}{Theorem}
\newtheorem{theorem}{Theorem}[section]

\newtheorem*{conjecture*}{Conjecture}
\newtheorem{lemma}[theorem]{Lemma}

\newtheorem*{corollary*}{Corollary}
\newtheorem{proposition}[theorem]{Proposition}

\advance \topmargin by -\headheight
\advance \topmargin by -\headsep
\evensidemargin \oddsidemargin
\begin{document}

\begin{centering}

{\huge \textbf{The structure of base phi expansions}}

\bigskip

{\bf \large F.~Michel Dekking}

\bigskip

{CWI, Amsterdam, and DIAM,  Delft University of Technology, Faculty EEMCS,\\ P.O.~Box 5031, 2600 GA Delft, The Netherlands.}

\medskip

{\footnotesize \it Email:  F.M.Dekking@TUDelft.nl}

\end{centering}

\medskip

\begin{abstract}
 \noindent In the base phi expansion any natural number is written uniquely as a sum of powers of the golden mean with coefficients 0 and 1, where it is required that the product of two consecutive digits is always 0. We tackle the problem of describing how these expansions look like. We classify the positive parts of the base phi expansions according to their suffices, and the negative parts according to their prefixes, specifying the sequences of occurrences of these digit blocks. Here the situation is much more complex than for the  Zeckendorf expansions, where any natural number is written uniquely as a sum of Fibonacci numbers with coefficients 0 and 1, where, again, it is required that the product of two consecutive digits is always 0. In a previous work we have classified the Zeckendorf expansions according to their suffices. It turned out that if we consider the suffices as labels on the Fibonacci tree, then the  numbers with a given suffix in their Zeckendorf expansion appear as generalized Beatty sequences in a natural way on this tree.
 We prove that the positive parts of the base phi expansions are a subsequence of the sequence of Zeckendorf expansions, giving an explicit formula in terms of a generalized Beatty sequence. The negative parts of the base phi expansions no longer appear lexicographically. We prove that all allowed digit blocks appear, and determine the order in which they do appear.
\end{abstract}

\medskip


\bigskip

\section{Introduction}

Let the golden mean be given by  $\varphi:=(1+\sqrt{5})/2$.\\
  Ignoring leading and trailing zeros, any  natural number $N$ can be written uniquely as
  $$N= \sum_{i=-\infty}^{\infty} d_i \varphi^i,$$
  with digits $d_i=0$ or 1, and where $d_id_{i+1} = 11$ is not allowed.
  As usual, we denote the base phi expansion of $N$ as $\beta(N)$, and we  write these expansions with a radix point as
  $$\beta(N) = d_{L}d_{L-1}\dots d_1d_0\cdot d_{-1}d_{-2} \dots d_{R+1}d_R.$$

 We define $$\beta^+(N)=d_{L}d_{L-1}\dots d_1d_0\; {\;\rm and\;}\; \beta^-(N)=d_{-1}d_{-2} \dots d_{R+1}d_R.$$
 So  $\beta(N)=\beta^+(N)\cdot\beta^-(N)$. For example, $\beta(2)=10\cdot 01$, and $\beta(3)=100\cdot 01$.

 \medskip

This paper deals with the following question: what are the words of 0's and 1's that can occur as digit blocks in the base phi expansion $N$, and for which numbers $N$ do they occur?

\medskip

In Section \ref{sec:phi}, we perform this task for the suffices of the $\beta^+$-part of the base phi expansions, and in Section \ref{sec:neg} for the complete $\beta^-\!$-part of the base phi expansions, and the prefixes of the $\beta^-\!$-part  of length at most 3.

\medskip

In Section \ref{sec:embed}, we establish in Theorem \ref{th:Zeckphi} a relationship between the  base phi expansions and Zeckendorf expansions, also known as Fibonacci representations.
This will permit us to exploit the results of the  paper \cite{Dekk-Zeck-structure} in Section \ref{sec:phi}.  See the paper \cite{Frougny-Saka} for a less direct approach, in terms of two-tape automata.

In Section \ref{sec:RST} we recall the recursive structure of base phi expansions, and derive some tools from this which will be useful in the final two sections.

In Section \ref{sec:closer} we take a closer look at the Lucas intervals.

In Section \ref{sec:GBS} we introduce generalized Beatty sequences, which for the base phi expansion take over the role played by arithmetic sequences in the classical expansions in base $b$, where $b$ is an integer larger than 1.

\medskip

We end this introduction by pointing out that there is a neat way to obtain $N$ from the $\beta^+(N)$-part of $\beta(N)$, without knowing the $\beta^-(N)$-part.
If $\beta(N)=\beta^+(N)\cdot\beta^-(N)$ is the base phi expansion of a natural  number $N$, then $N=\lceil \beta^+(N)\rceil$. Here $\lceil \cdot\rceil$ is the ceiling function.

For a proof, add the maximum number of powers corresponding to $\beta^-(N)$, taking into account that no 11 appears.
This is bounded by the geometric series starting at $\varphi^{-1}$ with common ratio $\varphi^{-2}$, i.e., by $\varphi^{-1}/(1-\varphi^{-2})=1.$

\section{Embedding base phi into Zeckendorf}\label{sec:embed}

We define the Lekkerkerker-Zeckendorf expansion. Let $(F_n)$ be the Fibonacci numbers. Let $\ddot{F}_0=1, \ddot{F}_1=2, \ddot{F}_2=3,\dots$ be the twice shifted Fibonacci numbers, defined by $\ddot{F}_i=F_{i+2}$.
  Ignoring leading and trailing zeros, any  natural number $N$ can  be written uniquely as
  $$N= \sum_{i=0}^{\infty} e_i \ddot{F}_i,$$
  with digits $e_i=0$ or 1, and where $e_ie_{i+1} = 11$ is not allowed.
  We denote the Zeckendorf expansion of $N$ as $Z(N)$.

  \medskip

  Let $V$ be the generalized Beatty sequence (cf. \cite{GBS}) defined by
  $$ V(n) = 3\lfloor n\varphi \rfloor + n +1.$$

\noindent   Here $\lfloor \cdot \rfloor$ denotes the floor function,  and $(\lfloor n\varphi\rfloor)$ is the well known lower Wythoff sequence.

  We define the function $S$ by
  $$S(n)=\max\{k\in \mathbb{N}: V(k)\le n\}-1.$$

\begin{theorem}\label{th:Zeckphi}
\noindent For all $N\ge 0$
$$\beta^+(N)=Z(N+S(N)).$$
\end{theorem}

\medskip

\noindent This theorem will be proved in the  Section \ref{sec:Pf-Zeckphi}.

\medskip

The basis for the embedding of the $\beta^+(N)$ into the collection of Zeckendorf words is the following analysis.

\subsection{The art of adding 1}\label{sec:add}

It is essential to give ourselves the freedom to write also non-admissible expansions in the form
$$\beta(N) = d_{L}d_{L-1}\dots d_1d_0\cdot d_{-1}d_{-2} \dots d_{R+1}d_R.$$
 For example, since $\beta(4) =101.01$ and
$\beta(2)=10\cdot 01$, we can write
\begin{equation}\label{eq:plus1} \beta(5)\doteq \beta(4)+1\doteq 101\cdot01+1\cdot 0\doteq 102\cdot01\doteq 110\cdot02\doteq1000\cdot1001.
\end{equation}
 Here the $\doteq$-sign indicates that we consider a non-admissible expansion.

It is convenient to generate all Zeckendorf expansions and base phi expansions by repeatedly adding the number 1.

When we compute $\beta(N)+1$ for some number $N$, then, in general, there is a carry both to the left and (two places) to the right.
This is illustrated by the example in Equation (\ref{eq:plus1}).
Note that there is not only a {\it double carry}, but that we also have to get rid of the 11's, by replacing them with 100's.
This is allowed because of the equation $\varphi^{n+2}=  \varphi^{n+1}+\varphi^{n}.$ We call this operation a {\it golden mean shift}.

When we compute $Z(N)+1$ for some number $N$, then we have to distinguish between $e_0=0$ and $e_0=1$:
$$Z(N)=e_L\dots e_2e_1\,0 {\quad\rm gives\quad} Z(N)+1=e_L\dots e_2e_1\,1$$
and
$$Z(N)=e_L\dots e_2e_1\,1 {\quad\rm gives\quad} Z(N)+1\doteq e_L\dots e_2\,10.$$
Here we used the $\doteq$-sign because (several) golden mean shifts might follow, where for the Zeckendorf expansion these are justified by the equation $F_{n+2}=F_{n+1}+F_n$. Note that replacing $e_11+1$ by $10$ follows from 1+1=2 (!).

\medskip

\noindent  For the convenience of the reader we provide a list of the Zeckendorf and base phi expansions of the first 18 natural numbers:

\bigskip

 \begin{tabular}{|c|c|c|}
   \hline
  \; $N^{\phantom{|}}$ & $Z(N)$ & $\beta(N)$  \\[.0cm]
   \hline
   1\; & \;\;\;\;\;\;\;\,1   & \;${1}\cdot$              \\
   2\; & \;\;\;\;\;\;10   & \;\:\,\,$1{0}\cdot01$     \\
   3\; &\;\;\; 100   & \;$10{0}\cdot01$          \\
   4\; &\;\;\; 101   & \;$10{1}\cdot01$          \\
   5\; & \;\;1000   & \;\:\,$100{0}\cdot1001$   \\
   6\; & \;\;1001   & \;\:\,$101{0}\cdot0001$   \\
   7\; & \;\;1010   & \;$1000{0}\cdot0001$      \\
   8\; & 10000    & \;$1000{1}\cdot 0001$     \\
   9\; & 10001    & \;$1001{0}\cdot0101$    \\
   \hline
 \end{tabular}\qquad
 \begin{tabular}{|c|c|c|}
   \hline
  \; $N^{\phantom{|}}$ & $Z(N)$ & $\beta(N)$  \\[.0cm]
   \hline
   10\; & \;\:10010  & \;$1010{0}\cdot0101$   \\
   11\; & \;\:10100 & \;$1010{1}\cdot0101$   \\
   12\; & \;\:10101 & \;\,\,$10000{0}\cdot101001$      \\
   13\; & 100000  & \;\,\,$10001{0}\cdot001001$      \\
   14\; & 100001  & \;\,\,$10010{0}\cdot001001$      \\
   15\; & 100010  & \;\,\,$10010{1}\cdot001001$   \\
   16\; & 100100  &  \;\, $10100{0}\cdot100001$\\
   17\; & 100101    &  \;\, $101010\cdot000001$\\
   18\; & 101000     &  \;\, $100000{0}\cdot000001$\\
   \hline
 \end{tabular}\quad

\subsection{Proof of  Theorem \ref{th:Zeckphi}}\label{sec:Pf-Zeckphi}

 The essential ingredient of the proof is the following result from  \cite{Dekk-phi-FQ}, Theorem 5.1 and Remark 5.4. An alternative, short proof of the first part could be given with the Propagation Principle from Section \ref{sec:RST}.

 \begin{proposition}\label{prop:D-numbers} Let $\beta(N)=(d_i(N))$ be the base phi expansion of a natural number $N$. Then:\\[-.3cm]
 \hspace*{1.5cm} $d_1d_0\cdot d_{-1}(N)=10\cdot1$ never occurs,\\[.1cm]
 \hspace*{1.5cm}  $d_1d_0\cdot d_{-1}(N)=00\cdot1$ if and only if $N=3\lfloor n\varphi\rfloor + n + 1$ for some natural number $n$.
\end{proposition}

\noindent {\it Proof of Theorem \ref{th:Zeckphi}: } One observes that there are many $\beta(N)$'s such that $\beta^+(N)=Z(N')$ for some $N'$. Moreover, if this is the case, then also $\beta^+(N+1)=Z(N'+1)$, {\it except} if $d_{-1}(N)=1$ in  $\beta(N)$. Indeed, as long as $d_{-1}(N)=0$ adding 1 gives the same result for both the Zeckendorf and the positive part of the base phi expansion, as seen in the previous section. However, suppose
$$Z(N')=\beta^+(N), {\; \rm and\;} d_{-1}(N)=1.$$
Then, by Proposition \ref{prop:D-numbers},  $d_1d_0\cdot d_{-1}(N)=00\cdot1$,
and adding 1 to $N$ gives the expansion $\beta(N+1)$ with digit block $d_1d_0\cdot d_{-1}(N+1)=10\cdot0$. So $\beta^+(N+1)$ ends in exactly the same two digits as $Z(N'+2)$, and in fact $\beta^+(N+1)=Z(N'+2)$. This means that one Zeckendorf expansion has been skipped: that of $N'+1$. Every time a $d_{-1}(N)=1$ occurs, this skipping takes place. Since $Z(0)=\beta^+(0),\dots, Z(5)=\beta^+(5)$, this gives the formula $\beta^+(N)=Z(N+S(N))$, with
$S(n)=\max\{k\in \mathbb{N}: 3\lfloor k\varphi\rfloor + k \le n\}$, by the second statement of Proposition \ref{prop:D-numbers}. \hfill $\Box$

\section{The recursive structure of base phi expansions}\label{sec:RST}

The Lucas numbers $(L_n)=(2, 1, 3, 4, 7, 11, 18, 29, 47, 76,123, 199, 322,\dots)$ are defined by
$$  L_0 = 2,\quad L_1 = 1,\quad L_n = L_{n-1} + L_{n-2}\quad {\rm for \:}n\ge 2.$$
The Lucas numbers have a particularly simple base phi expansion.

\noindent From  the well-known formula
$L_{2n}=\varphi^{2n}+\varphi^{-2n}$, and the recursion $L_{2n+1}=L_{2n}+L_{2n-1}$ we have for all $n\ge 1$
\begin{equation}\label{eq:Lm}
 \beta(L_{2n}) = 10^{2n}\cdot0^{2n-1}1,\quad \beta(L_{2n+1}) = 1(01)^n\cdot(01)^n.
\end{equation}
By iterated application of the double carry and the golden mean shift  to $\beta(L_{2n+1})+\beta(1)$,\; and a similar operation for  $\beta(L_{2n+2}-1)$   (see also the last page of \cite{Dekk-How-to-add}) one finds that for all $n\ge 1$

\begin{equation}\label{eq:Lmplus1}
 \beta(L_{2n+1}+1) = 10^{2n+1}\cdot(10)^n01,\quad  \beta(L_{2n+2}-1)=(10)^{n+1}\cdot 0^{2n+1}1.
\end{equation}

\noindent As in \cite{Dekk-phi-FQ} we partition the natural numbers into Lucas intervals\:
$$\Lambda_{2n}:=[L_{2n},\,L_{2n+1}] \quad{\rm and\quad} \Lambda_{2n+1}:=[L_{2n+1}+1,\, L_{2n+2}-1].$$
The basic idea behind this partition is that if
 $$\beta(N) = d_{L}d_{L-1}\dots d_1d_0\cdot d_{-1}d_{-2} \dots d_{R+1}d_R,$$
then the left most index $L=L(N)$ and the right most index $R=R(N)$ satisfy
$$L(N)=|R(N)|=2n \;{\rm iff}\; N\in \Lambda_{2n}, \quad L(N)=2n\!+1,\; |R(N)|=2n\!+2 \;{\rm iff}\; N\in \Lambda_{2n+1}.$$
This is not hard to see from the simple expressions we have for the $\beta$-expansions of the Lucas numbers; see also Theorem 1 in \cite{Grabner94}.

\medskip

To obtain recursive relations, the interval $\Lambda_{2n+1}=[L_{2n+1}+1, L_{2n+2}-1]$ has to be divided into three subintervals. These three intervals are\\[-.8cm]
 \begin{align*}
I_n:=&[L_{2n+1}+1,\, L_{2n+1}+L_{2n-2}-1],\\
J_n:=&[L_{2n+1}+L_{2n-2},\, L_{2n+1}+L_{2n-1}],\\
K_n:=&[L_{2n+1}+L_{2n-1}+1,\, L_{2n+2}-1].
\end{align*}

\noindent It will be very convenient to use the free group versions of words of 0's and 1's. So, for example, $(01)^{-1}0001=1^{-1}001$.

\begin{theorem}{\bf [Recursive structure theorem]}\label{th:rec}

\noindent{\,\bf I\;} For all $n\ge 1$ and $k=0,\dots,L_{2n-1}$
one has $ \beta(L_{2n}+k) =  \beta(L_{2n})+ \beta(k) = 10\dots0 \,\beta(k)\, 0\dots 01.$
\noindent{\bf II} For all $n\ge 2$ and $k=1,\dots,L_{2n-2}-1$
\begin{align*}
I_n:&\quad \beta(L_{2n+1}+k) = 1000(10)^{-1}\beta(L_{2n-1}+k)(01)^{-1}1001,\\ K_n:&\quad\beta(L_{2n+1}+L_{2n-1}+k)=1010(10)^{-1}\beta(L_{2n-1}+k)(01)^{-1}0001.
\end{align*}
Moreover, for all $n\ge 2$ and $k=0,\dots,L_{2n-3}$
$$\hspace*{0.7cm}J_n:\quad\beta(L_{2n+1}+L_{2n-2}+k) = 10010(10)^{-1}\beta(L_{2n-2}+k)(01)^{-1}001001.$$
\end{theorem}

\medskip

See \cite{Dekk-How-to-add} for a proof of this theorem.

\medskip

As an illustration of the use of Theorem \ref{th:rec} we shall now prove a lemma that we need in Section \ref{sec:phi}.

\begin{lemma}\label{lem:no} Let $m\ge 1$ be an integer. There are {\bf (a)} no expansions $\beta(N)$ with the digit block $d_{2m}\dots d_0\cdot d_{-1}(N)=10^{2m}\cdot1$, and  there are {\bf (b)} no expansions $\beta(N)$ with the digit block   $d_{2m+1}\dots d_0\cdot d_{-1}(N)=10^{2m+1}\cdot0$.
\end{lemma}

\noindent{\it Proof:} \:{\bf (a)}. The first time $d_{2m}\dots d_0=10^{2m}$ occurs is for $N=L_{2m}$, and then $d_{-1}(N)=0$ (see $\beta(L_{2m})$ formula above). This is also the only occurrence of the digit block $10^{2m}$ at the end of the expansions of the numbers $N$ in $\Lambda_{2m}$. It is also obvious that the digit block $10^{2m}$ will not appear at the end of the expansions of the numbers $N$ in $\Lambda_{2m+1}$.

From part {\bf I} of the Recursive Structure Theorem we see that the digit block $10^{2m}$ at the end of the expansions of the numbers $N$ in $\Lambda_{2m+2}$ only occurs in combination with $d_{-1}(N)=0$.

From part {\bf II} of the Recursive Structure Theorem we will see that the digit block $10^{2m}$ at the end of the expansions of the numbers $N$ in $\Lambda_{2m+3}$ only occurs in combination with $d_{-1}(N)=0$. This is definitely more complicated than this observation for $\Lambda_{2m+2}$. We have to split $\Lambda_{2m+3}$ into the three pieces $I_{m+1}, J_{m+1}$ and $K_{m+1}$. The middle piece $J_{m+1}$ corresponds to numbers in $\Lambda_{2m}$, from which we  already know that $d_{2m}\dots d_0(N)=10^{2m}$ implies that $d_{-1}(N)=0$. The numbers $N$ in the first piece, $I_{m+1}$, correspond to numbers in $\Lambda_{2m+1}$ from which the digits $d_{2m+1}d_{2m}=10$ have been replaced by
the digits $d_{2m+3}d_{2m+2}d_{2m+1}d_{2m}=1000$. In particular $d_{2m}=0$ excludes any occurrence of $d_{2m}\dots d_0=10^{2m}$. In the same way occurrences of $d_{2m}\dots d_0=10^{2m}$ in $K_{m+1}$ are excluded.

The final conclusion is that both intervals $\Lambda_{2m+2}$ and $\Lambda_{2m+3}$ only contain numbers $N$ for which the occurrence of  $10^{2m}$ as end block implies $d_{-1}(N)=0$. In the same way, these properties of $\Lambda_{2m+2}$ and $\Lambda_{2m+3}$ carry over to the two Lucas intervals $\Lambda_{2m+4}$ and $\Lambda_{2m+5}$, and we can finish the proof by induction.

{\bf (b)}. The first time $d_{2m}\dots d_0=10^{2m+1}$ occurs is for $N=L_{2m+1}+1$ in $\Lambda_{2m+1}$ , and then $d_{-1}(N)=1$ (see Equation (\ref{eq:Lmplus1})). This is also the only occurrence of $d_{2m}\dots d_0=10^{2m+1}$ in $\Lambda_{2m+1}$. Moreover, in $\Lambda_{2m+2}$  the word $10^{2m+1}$  does not occur at all as end block. We finish the proof as in Part {\bf (a)}, with the sole difference that now $1010^{2m+1}$ occurring as end block in $\Lambda_{2m+3}$,
yields an instance of $10^{2m+1}\cdot 1$ in $\Lambda_{2m+3}$. \hfill $\Box$

\bigskip

It is convenient to have a second version of the Recursive Structure Theorem which involves a higher resemblance between the even Part {\,\bf I\;} case, and the odd Part {\,\bf II\;}. It will also be convenient to have the $\Lambda$-intervals play a more visible role in the recursion. In fact, it is easy to check that the three intervals $I_n, J_n$ and $K_n$ in the  Recursive Structure Theorem satisfy
$$I_n=\Lambda^{(a)}_{2n-1}:=\Lambda_{2n-1}+L_{2n},\; J_n=\Lambda^{(b)}_{2n-2}:=\Lambda_{2n-2}+L_{2n+1},\; K_n=\Lambda^{(c)}_{2n-1}:=\Lambda_{2n-1}+L_{2n+1}. $$
In this equation we employ the usual notation $A+x:=\{a+x:a\in A\}$ for a set of real numbers $A$ and a real number $x$.

\begin{theorem}{\bf [Recursive structure theorem: 2nd version]}\label{th:rec2}\\
\noindent{\,\bf (i):  Odd\;} For all $n\ge 1$ one has\\[-.3cm]
$$\Lambda_{2n+1}=\Lambda^{(a)}_{2n-1}\cup\Lambda^{(b)}_{2n-2}\cup\Lambda^{(c)}_{2n-1}, $$
where $\Lambda^{(a)}_{2n-1}=\Lambda_{2n-1}+L_{2n}$,\; $\Lambda^{(b)}_{2n-2}=\Lambda_{2n-2}+L_{2n+1}$, and $\Lambda^{(c)}_{2n-1}=\Lambda_{2n-1}+L_{2n+1}$.\\
We have\\[-.8cm]
\begin{subequations} \label{eq:shift-odd}
\begin{align}
  \beta(N)= & \;1000(10)^{-1}\,\beta(N-L_{2n})\,(01)^{-1}1001& for\; N\in \Lambda^{(a)}_{2n-1}, \\
  \beta(N)= & \; 100\,\beta(N-L_{2n+1})\,(01)^{-1}001001&      for\; N\in  \Lambda^{(b)}_{2n-2},\\
  \beta(N)= & \;10\,\beta(N-L_{2n+1})\,(01)^{-1}0001 &         for\; N\in \Lambda^{(c)}_{2n-1}.
\end{align}\\[-.6cm]
\end{subequations}
\noindent{\,\bf (ii): Even\;} For all $n\ge 1$ one has\\[-.3cm]
 $$\Lambda_{2n+2}=\Lambda^{(a)}_{2n}\cup\Lambda^{(b)}_{2n-1}\cup\Lambda^{(c)}_{2n}, $$
 where $\Lambda^{(a)}_{2n}=\Lambda_{2n}+L_{2n+1}$,\; $\Lambda^{(b)}_{2n-1}=\Lambda_{2n-1}+L_{2n+2}$, and $\Lambda^{(c)}_{2n}=\Lambda_{2n}+L_{2n+2}$.\\
 We have\\[-.8cm]
\begin{subequations} \label{eq:shift-even}
\begin{align}
  \beta(N)= & \;1000(10)^{-1}\,\beta(N-L_{2n+1})\,(01)^{-1}0001 & for\; N\in \Lambda^{(a)}_{2n},\phantom{x} \label{eq:5a} \\
  \beta(N)= & \; 100\,\beta(N-L_{2n+2})\,01 &                     for\; N\in \Lambda^{(b)}_{2n-1},  \label{eq:5b}\\
  \beta(N)= & \;10\,\beta(N-L_{2n+1})\,01 &                       for\; N\in \Lambda^{(c)}_{2n}.\phantom{x.}  \label{eq:5c}
\end{align}
\end{subequations}
\end{theorem}

\noindent{\it Proof:} \:{\bf (i):  Odd\;} This is a rephrasing of {Part (II) in Theorem \ref{th:rec}.\\
{\bf (ii):  Even\;} We  start by showing that the three intervals $\Lambda^{(a)}_{2n},\Lambda^{(b)}_{2n-1},\Lambda^{(c)}_{2n}$ partition $\Lambda_{2n+2}$.

The first number in $\Lambda^{(a)}_{2n}$ is  $L_{2n}+L_{2n+1}=L_{2n+2}$, which is the first number of $\Lambda_{2n+2}$. The last number in $\Lambda^{(a)}_{2n}$ is  $L_{2n+1}+L_{2n+1}=2L_{2n+1}$.

The first number in $\Lambda^{(b)}_{2n-1}$ is  $L_{2n-1}+1+L_{2n+2}=L_{2n-1}+1+L_{2n}+L_{2n+1}=2L_{2n+1}+1$, which indeed, is the successor of the last number in $\Lambda^{(a)}_{2n}$.

The last number in $\Lambda^{(b)}_{2n-1}$ is  $L_{2n}-1+L_{2n+2}$, which indeed has successor $L_{2n}+L_{2n+2}$, the first number in $\Lambda^{(c)}_{2n}$. Finally, the last number in $\Lambda^{(c)}_{2n}$ is $L_{2n+1}+L_{2n+2}=L_{2n+3}$, which is the last number in $\Lambda_{2n+2}$.

\medskip

To prove Equation (\ref{eq:5a}), we first show, using Equation (\ref{eq:Lm}) twice, that this equation is correct for $N=L_{2n+2}$, which is the first number of $\Lambda^{(a)}_{2n}$:
\begin{align*}
\beta(L_{2n+2}) & = 10^{2n+2}\cdot0^{2n+1}1\\
                & = 1000\,0^{2n-1}\cdot0^{2n-2}\,0001 \\
                & = 1000(10)^{-1}10^{2n}\cdot0^{2n-1}1\,(01)^{-1}0001 \\
                & = 1000(10)^{-1}\beta(L_{2n})\,(01)^{-1}0001 \\
                & =  1000(10)^{-1}\beta(L_{2n+2}-L_{2n+1})\,(01)^{-1}0001.
\end{align*}

Equation (\ref{eq:5a}) will also be correct for all other $N\in\Lambda^{(a)}_{2n}$, because as above, the digit block $d_Ld_{L-1}d_{L-2}d_{L-3}(N)$ will always be 1000,
and the digit block $d_{L-2}d_{L-3}(N-L_{2n+1})$ will always be 10. For the negative digits we have a similar property.

Equation (\ref{eq:5b}) follows directly from the fact that if $N\in \Lambda^{(b)}_{2n-1}$, then
\begin{align*}
\beta(N-L_{2n+2})+\beta(L_{2n+2}) & = d_{2n-1}\dots d_0\cdot d_{-1}\dots d_{-2n} + 10^{2n+2}\cdot0^{2n+1}1\\
                                  & = d_{2n-1}\dots d_0\cdot d_{-1}\dots d_{-2n} + 100\,0^{2n}\cdot0^{2n}\,01\\
                                  & = 100d_{2n-1}\dots d_0\cdot d_{-1}\dots d_{-2n}01,
\end{align*}
since the numbers in $\Lambda_{2n-1}$ have a $\beta$-expansion $d_{2n-1}\dots d_0\cdot d_{-1}\dots d_{-2n}$ with $2n$ digits on the left and $2n$ digits on the right.
Note that we do not have to use the $\dot{=}$-sign as there are no double carries or golden mean shifts.

Equation (\ref{eq:5c}) follows in the same way. \hfill $\Box$

\bigskip

Lemma \ref{lem:no} is an example of a general phenomenon, which we call the Propagation Principle. It has an extension to combinations of digit blocks which we will give in Lemma \ref{lem:prop}.

The Propagation Principle is closely connected to the following notion.
We say an interval $\Gamma$ and a union of intervals  $\Delta$  of natural numbers are \emph{$\beta$-congruent modulo $q$} for some natural number $q$ if
$\Delta$ is a disjoint union of translations of $\Lambda$-intervals, such that for all $j=1,\dots,|\Gamma|$,
 if $N$ is the $j^{\rm th}$ element of  $\Gamma$, and $N'$ is the $j^{\rm th}$ element of $\Delta$, then
$$d_{q-1}\dots d_1d_0\cdot d_{-1}\dots d_{-q}(N)=d_{q-1}\dots d_1d_0\cdot d_{-1}\dots d_{-q}(N').$$
We write this as\,   $\Gamma\cong \Delta_1\Delta_2\dots\Delta_r \mod q$ when the number of translations of $\Lambda$-intervals in $\Delta$ equals $r$.
Note that the definition implies that the $r$ disjoint translations of $\Lambda$-intervals appear in the natural order, and that we refrain from indicating the translations.

Simple examples are $\Lambda_5 \cong \Lambda_3\Lambda_2\Lambda_3 \mod 1$  and $\Lambda_6 \cong \Lambda_4\Lambda_3\Lambda_4 \mod 3$.
 Theorem \ref{th:rec2} is a source of many more examples.

\medskip

An important observation is that if $\Gamma\cong \Delta_1\Delta_2\dots\Delta_r \mod q$ and $\Gamma'\cong \Delta'_1\Delta'_2\dots\Delta'_{r'} \mod q'$, and $\Gamma\cup \Gamma'$ is an interval, then
\begin{equation}\label{eq:GG}
 \Gamma\Gamma':=\Gamma\cup \Gamma'\cong \Delta_1\Delta_2\dots\Delta_r\Delta'_1\Delta'_2\dots\Delta'_{r'} \mod \min\{q,q'\}.
\end{equation}

\medskip

To keep the formulation and the proof of the following lemma simple, we only formulate it for central digit blocks of length 8 (i.e., $q=4$).
In the following, occurrences of digit blocks in $\beta$-expansions have to be interpreted with additional $0$'s added to the left of the expansion.

\begin{lemma}{\bf [Propagation Principle]}\label{lem:prop}

\noindent {\bf (a)}\, Suppose the digit block $d_{3}\dots d_0\cdot d_{-1}\dots d_{-4}$, does not occur in the $\beta$-expansions of the numbers $N=1,2,\dots,17$. Then it does not occur in any $\beta$-expansion.

\noindent {\bf (b)}\, Let $D$ be an integer between 1 and 4. Suppose the digit block $d_{3}\dots d_0\cdot d_{-1}\dots d_{-4}$ occurs in the $\beta$-expansion of $N$ if and only if the digit block $e_{3}\dots e_0\cdot e_{-1}\dots e_{-4}$ occurs in the $\beta$-expansion of the number $N-D$, for $N=D,D+1,\dots,D+17$.    Then this coupled occurrence holds for all $N$.
\end{lemma}

\noindent{\it Proof:} \:{\bf (a)}\, Let us say that a Lucas interval $\Lambda_m$ satisfies property $\cal D$ if the digit block $d_{3}\dots d_0\cdot d_{-1}\dots d_{-4}$, does not occur in the $\beta$-expansions  of the numbers $N$ from $\Lambda_m$. Note that $N=17$ is the last number in  $\Lambda_5$, so it is given that the intervals $\Lambda_1,\dots,\Lambda_5$ all satisfy property $\cal D$. Also $\Lambda_6$ satisfies property $\cal D$, by an application of Theorem \ref{th:rec}, Part {\bf I}.

The interval $\Lambda_7= \Lambda^{(a)}_{5}\cup\Lambda^{(b)}_{4}\cup\Lambda^{(c)}_{5}$ satisfies property $\cal D$. For $\Lambda^{(a)}_{5}$, this follows since $\Lambda_5$ satisfies property $\cal D$, and (\ref{eq:5a}) does not change the central block of length 8. The same argument applies to $\Lambda^{(c)}_{5}$.
For the interval $\Lambda^{(b)}_{4}$,  Equation (\ref{eq:5b}) gives that the positive digit blocks $d_{3}\dots d_0$ are the same as for the corresponding numbers in $\Lambda_4$, and that  the negative digit blocks are $d_{-1}\dots d_{-4}(7)(01)^{-1}00=0000$ and $d_{-1}\dots d_{-4}(9)(01)^{-1}00=0100$, which already occurred in the expansions $\beta(0)$ and $\beta(3)$.

The interval $\Lambda_8= \Lambda^{(a)}_{6}\cup\Lambda^{(b)}_{5}\cup\Lambda^{(c)}_{6}$ also satisfies property $\cal D$, since the word transformations in Equation (\ref{eq:shift-even})  do not change the central blocks of length 8 in $\Lambda_{6}$, nor in $\Lambda_{5}$.
Another way to put this, is that $\Lambda_8\cong \Lambda_{6}\Lambda_{5}\Lambda_{6} \mod 4$. Since the $\beta$-expansions only get longer, we have in fact that
$\Lambda_m\cong \Lambda_{m-2}\Lambda_{m-3}\Lambda_{m-2} \mod 4$ for all $m\ge 8$. Thus it follows by induction that $\Lambda_m$ satisfies property $\cal D$ for all $m\ge 8$.

\noindent\:{\bf (b)}\, Let us say that a Lucas interval $\Lambda_m$, $m\ge 1$  satisfies property $\cal E$ if the numbers $N$ from $\Lambda_m$ have the property that the digit block $d_{3}\dots d_0\cdot d_{-1}\dots d_{-4}$ occurs in the $\beta$-expansion of $N$ if and only if the digit block $e_{3}\dots e_0\cdot e_{-1}\dots e_{-4}$ occurs in the $\beta$-expansion of $N-D$. Then it is given that $\Lambda_1,\dots,\Lambda_5$ all satisfy property $\cal E$. The proof continues as in part {\bf (a)}, but we have to take into account that the numbers $N-D$ and $N$ can be elements of different Lucas intervals. This `boundary' problem is easily solved by induction: it is given for $\Lambda_4\Lambda_{5}$  and $\Lambda_5\Lambda_{6}$, and the equation used for induction is\\[-.4cm]
$$\Lambda_{m+1}\Lambda_{m+2}\cong \Lambda_{m-1}\Lambda_{m-2}\Lambda_{m-1}\Lambda_{m}\Lambda_{m-1}\Lambda_{m} \mod 4.$$
This equation is an instance of Equation (\ref{eq:GG}).\hfill $\Box$

\section{A closer look at the Lucas intervals}\label{sec:closer}

Here we say more on the idea of  splitting Lucas intervals in unions of translated Lucas intervals.
To keep the presentation simple, we start with showing how all the natural numbers can be split into translations of the three Lucas intervals $\Lambda_3, \Lambda_4$ and $\Lambda_5$.

This can of course be done in many ways, but we will consider a way derived from the Recursive Structure Theorem \ref{th:rec2}.
One has
\begin{align*}
  \Lambda_6= & \:\Lambda^{(a)}_{4}\cup\Lambda^{(b)}_{3}\cup\Lambda^{(c)}_{4}=[\Lambda_{4}\!+\!L_5]\cup[\Lambda_{3}\!+\!L_6]\cup[\Lambda_{4}\!+\!L_6], \\
  \Lambda_7= &\: \Lambda^{(a)}_{5}\cup\Lambda^{(b)}_{4}\cup\Lambda^{(c)}_{5}=[\Lambda_{5}\!+\!L_5]\cup[\Lambda_{4}\!+\!L_7]\cup[\Lambda_{5}\!+\!L_7],  \\
  \Lambda_8= & \:\Lambda^{(a)}_{6}\cup\Lambda^{(b)}_{5}\cup\Lambda^{(c)}_{6}=[\Lambda_{6}\!+\!L_7]\cup[\Lambda_{5}\!+\!L_8]\cup[\Lambda_{6}\!+\!L_8]\\
  = & \: [\Lambda_{4}\!+\!L_5\!+\!L_7]\cup[\Lambda_{3}\!+\!L_6\!+\!L_7]\cup[\Lambda_{4}\!+\!L_6\!+\!L_7]\cup[\Lambda_{5}\!+\!L_8]\\
     & \hspace*{6cm}\cup[\Lambda_{4}\!+\!L_5\!+\!L_8]\cup[\Lambda_{3}\!+\!L_6\!+\!L_8]\cup[\Lambda_{4}\!+\!L_6\!+\!L_8].
\end{align*}
Note how the splitting of $\Lambda_6$ was used in the splitting of $\Lambda_8$. Continuing in this fashion we obtain inductively a splitting of all Lucas intervals $\Lambda_n$, which we call the \emph{canonical} splitting.

What is the sequence of translated intervals $\Lambda_3, \Lambda_4$ and $\Lambda_5$ created in this way?

\medskip
Let the word $C(\Lambda_{n})$ code these successive intervals in $\Lambda_n$ by their indices 3, 4 or 5. Let $\kappa$ be the morphism on the monoid $\{3,4,5\}^*$ defined by
$$\kappa(3)=5, \quad \kappa(4)=434 \quad \kappa(5)=545.$$

\medskip

\begin{theorem}\label{th:L-345}  For any $n\ge 3$ the interval $\Lambda_n$ is a union of adjacent translations of the three intervals $\Lambda_3, \Lambda_4$ and $\Lambda_5$.  If $C(\cdot)$ is the coding function for the canonical splittings  then for $n\ge 0$
$$C(\Lambda_{2n+4})=\kappa^n(4), \quad C(\Lambda_{2n+5})=\kappa^n(5).$$
\end{theorem}

\noindent{\it Proof:} By induction. For $n=0$ this is trivially true.

Suppose it is true for $k =1, \dots n$. Then by Theorem \ref{th:rec2},
\begin{align*}
  C(\Lambda_{2n+6})&=C(\Lambda_{2n+4})C(\Lambda_{2n+3})C(\Lambda_{2n+4})=\kappa^n(4)\kappa^{n-1}(5)\kappa^n(4)=\kappa^{n-1}(\kappa(4)5\kappa(4))\\
  &=\kappa^{n-1}(4345434)=\kappa^{n-1}(\kappa^2(4))=\kappa^{n+1}(4), \\
  C(\Lambda_{2n+7}) &=C(\Lambda_{2n+5})C(\Lambda_{2n+4})C(\Lambda_{2n+5})= \kappa^n(5)\kappa^n(4)\kappa^n(5)=\kappa^{n}(545)=\kappa^{n+1}(5).\quad \Box
\end{align*}

\medskip

We continue this analysis, now focussing on the partition of the natural numbers by the intervals
$$\Xi_n:=\Lambda_{2n-1}\cup\Lambda_{2n}=[L_{2n-1}+1,L_{2n+1}].$$
The relevance of the $\Xi_n$ is that these are exactly the intervals where $\beta^-(N)$ has length $2n$, for $n \ge 1$. The results in the sequel of this section will therefore be useful in Section \ref{sec:neg}.

\medskip

There are three (Sturmian) morphisms $f,g$ and $h$ that play an important role in these results, where it is convenient to look at $a$ and $b$ both as integers and as abstract letters. The morphisms are given by
\begin{equation}\label{eq:Fib3}
f:\mor{aba}{ab}\,,\qquad g:\mor{baa}{ba}\,,\qquad h:\mor{aab}{ab}.
\end{equation}



\begin{theorem}\label{th:d-345}  For any $n\ge 2$ the interval $\Xi_n$ is a union of adjacent translations of the three intervals $\Lambda_3, \Lambda_4$ and $\Lambda_5$.  If $C(\cdot)$ is the coding function for the canonical splittings,  then for $n\ge 0$
$$C(\Xi_{n+2})=\delta(h^n(b)),$$
where $\delta$ is the decoration morphism given by $\delta(a)=54,\:\delta(b)=34$.
\end{theorem}

\noindent{\it Proof:} We first establish the commutation relation\;$\kappa\,\delta=\delta \,h.$\\
It suffices to prove this for the generators, and indeed:
$$  \kappa(\delta(a))=\kappa(54)=545434=\delta(aab)=\delta(h(a)),\quad  \kappa(\delta(b))=\kappa(34)=5434=\delta(ab)=\delta(h(b)). $$
Using Theorem \ref{th:L-345}, and the commutation relation we obtain for $n\ge 1$
\begin{align*}
  C(\Xi_{n+2})&=C(\Lambda_{2n+3})C(\Lambda_{2n+4})=\\
  &=\kappa^{n-1}(5)\kappa^{n}(4)=\kappa^{n-1}(5434)=\kappa^{n-1}(\delta(ab))=\delta(h^{n-1}(ab))=\delta(h^n(b)).
\end{align*}
For $n= 0$ we have $\Xi_2=\Lambda_3\cup\Lambda_4$, so $C(\Xi_2)=34=\delta(b)$. \quad $\Box$

\bigskip

\section{Generalized Beatty sequences}\label{sec:GBS}

Let $\alpha$ be an irrational number larger than 1. We call any sequence $V$ with terms of the form $V_n = p\lfloor n \alpha \rfloor + q n +r $, $n\ge 1$ a \emph{generalized Beatty sequence}. Here $p,q$ and $r$ are integers, called the \emph{parameters} of $V$, and we write $V=V(p,q,r)$.

In this paper we will only consider the case $\alpha=\varphi$, the golden mean, so any mention of a generalized Beatty sequence assumes that $\alpha=\varphi$.
A prominent role is played by the lower Wythoff sequence $A:=V(1,0,0)$ and the upper Wythoff sequence $B:=V(1,1,0)$. These are complementary sequences, associated to the Beatty pair ($\varphi,\varphi^2)$.

Here is the key lemma that tells us how generalized Beatty sequences behave under compositions. In its statement below, as Lemma \ref{lem:VA}, a typo in its source is corrected.

\begin{lemma}{\bf (\cite{GBS}, Corollary 2) }\label{lem:VA} Let $V$ be a generalized Beatty sequence with parameters $(p,q,r)$. Then $VA$ and $VB$ are generalized Beatty sequences with parameters $(p_{V\!A},q_{V\!A},r_{V\!A})=(p+q,p,r-p)$ and $(p_{V\!B},q_{V\!B},r_{V\!B})=(2p+q,p+q,r)$.
\end{lemma}

It will be useful later on to have a sort of converse of this lemma.

If $C$ and $D$ are two $\N$-valued sequences, then we denote by $C\sqcup D$ the sequence whose terms give the set $C(\N)\cup D(\N)$, in increasing order.

\begin{lemma}\label{lem:VAconv} Let $V=V(p,q,r)$ be a generalized Beatty sequence. Let $U$ and $W$ be two disjoint sequences with union $V=U\sqcup W$:
$$U(\N)\cap W(\N)=\emptyset,\quad U(\N)\cup W(\N)=V(\N).$$
Suppose $U$ is a generalized Beatty sequence with parameters $(p+q, p, r-p)$. Then $W$ is the generalized Beatty sequence with parameters $(2p+q,p+q,r)$.
\end{lemma}

\noindent{\it Proof:} According to Lemma \ref{lem:VA}, we have $U=VA$.  Since $A$ and $B$ are disjoint with union $\N$, we must have $W=VB$, and Lemma \ref{lem:VA} gives that $W$ is a generalized Beatty sequence with parameters $(2p+q,p+q,r)$.\hfill$\Box$

\medskip

Here is the key lemma to `recognize' a generalized Beatty sequence, taken from \cite{GBS}.
If $S$ is a sequence, we denote its  sequence of first order differences as $\Delta S$, i.e., $\Delta S$ is defined by
$$\Delta S(n) = S(n+1)-S(n), \quad {\rm for\;} n=1,2\dots.$$

\begin{lemma}\label{lem:diff}{ \rm \bf(\cite{GBS})} Let $V = (V_n)_{n \geq 1}$ be the generalized Beatty
sequence defined by $V_n = p\lfloor n \varphi \rfloor + q n +r$, and let $\Delta V$ be the
sequence of its first differences. Then $\Delta V$ is the Fibonacci word on the alphabet
$\{2p+q, p+q\}$. Conversely, if $x_{a,b}$ is the Fibonacci word on the alphabet
$\{a,b\}$,  then any $V$ with $\Delta V= x_{a,b}$ is a generalized Beatty sequence
$V=V(a-b,2b-a,r)$ for some integer $r$.
\end{lemma}

\section{The positive powers of the golden mean}\label{sec:phi}

For any digit block $w$  we will determine the  sequence $R_w$ of those numbers $N$ with digit block $w=d_{m-1}\dots d_0$ as suffix of $\beta^+(N)$. We sometimes call $w$  an \emph{end block} of $\beta^+(N)$. More generally, we are also interested in occurrence sequences of numbers $N$ with $d_{m-1}\dots d_0(N)=w$ and $d_{-1}\dots d_{-m'}(N)=v$. We denote these as $R_{w\cdot v}$.

For a couple of small values of $m,m'$, we have the following result from the paper \cite{Dekk-phi-FQ}, Theorem 5.1.

\begin{theorem}\label{th:d0d1} {\bf (\cite{Dekk-phi-FQ})} Let $\beta(N)=(d_i(N))$ be the base phi expansion of a natural number $N$. Then:\\[.1cm]
  $R_{1}=V_0(1,2,1)$,\quad   $R_{10}=V(1,2,-1)$,\quad   $R_{00\cdot 0}=V_0(1,2,0)$,\quad   $R_{00\cdot 1}=V(3,1,1)$.
\end{theorem}

Here it made sense to add $N=1$ to $V(1,2, 1)$, and $N=0$ to  $R_{00\cdot 0}$. We accomplished this  by adding the $n=0$ term to the generalized Beatty sequence $V$: we define $V_0$ by $$V_0(p,q,r) := (p\lfloor n \phi \rfloor + q n +r)_{n\ge 0}.$$

\medskip

The digit blocks $w=d_{m-1}\dots d_1\,0$ behave rather differently from digit blocks $w=d_{m-1}\dots d_1\,1$. We therefore analyse these cases separately, in Section \ref{sec:w0} and \ref{sec:w1} .

\subsection{Digit blocks $w=d_{m-1}\dots d_10$}\label{sec:w0}

  We order the digit blocks $w$ with $d_0=0$ in a Fibonacci tree. The first four levels of this tree are depicted below.

\bigskip

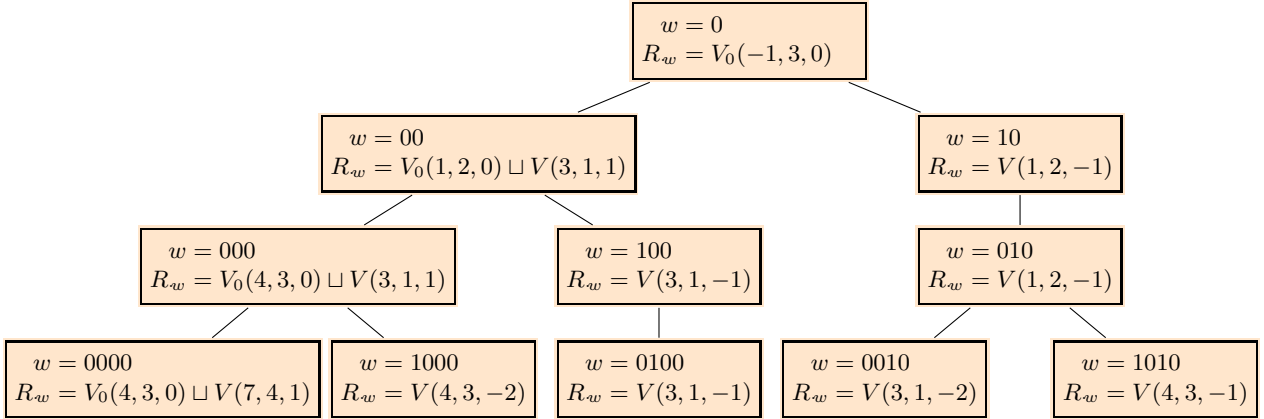
\hspace*{-1.1cm}
\begin{tikzpicture}
[level distance=15mm,
every node/.style={fill=orange!20,rectangle,inner sep=1pt},
level 1/.style={sibling distance=72mm,nodes={fill=orange!20}},
level 2/.style={sibling distance=48mm,nodes={fill=orange!20}},
level 3/.style={sibling distance=36mm,nodes={fill=orange!20}}]
\node {\footnotesize  $\duoZ{=0}{=V_0(-1,3,0)\quad}$}            
child {node {\footnotesize $\duoZ{=00}{=V_0(1,2,0) \sqcup V(3,1,1)}$}            
 child {node {\footnotesize $\duoZ{=000}{=V_0(4,3,0) \sqcup V(3,1,1)}$}          
  child {node {\footnotesize $\duoZ{=0000}{=V_0(4,3,0) \sqcup V(7,4,1)}$}}       
  child {node {\footnotesize $\duoZ{=1000}{=V(4,3,-2)}$}}     
}
 child {node {\footnotesize $\duoZ{=100}{=V(3,1,-1)}$}         
  child {node {\footnotesize $\duoZ{=0100}{=V(3,1,-1)}$}}     
}
}
child {node {\footnotesize $\duoZ{=10}{=V(1,2,-1)}$}            
 child {node {\footnotesize $\duoZ{=010}{=V(1,2,-1)}$}          
  child {node {\footnotesize $\duoZ{=0010}{=V(3,1,-2)}$}}     
   child {node {\footnotesize $\duoZ{=1010}{=V(4,3,-1)}$}}  
}
};
\end{tikzpicture}

\bigskip

We start with the short words $w$.

\medskip

\begin{proposition}\label{prop:small} The sequence of occurrences $R_w$ of numbers $N$ such that  the  digits $d_{m-1} \dots d_0$ of the base phi expansion of $N$ are equal to $w$, i.e., $d_{m-1}\dots d_0(N)=w$, is given for  the words $w$ of length at most 3, and ending in 0 by\\
\mbox{\rm a)} $R_0 = V(-1,3,0)$,\\
\mbox{\rm b)} $R_{00} = V_0(1,2,0) \sqcup V(3,1,1)$,\\
\mbox{\rm c)} $R_{10}=R_{010} = V(1,2,-1)$,\\
\mbox{\rm d)} $R_{000} = V_0(4,3,0) \sqcup V(3,1,1)$,\\
\mbox{\rm e)} $R_{100} = V(3,1,-1)$.
\end{proposition}

\noindent{\it Proof:} \noindent {\rm a)} $w=0$:  Since the numbers $\varphi+2$ and $3-\varphi$ form a Beatty pair, i.e.,
$$\frac{1}{\varphi+2}+\frac{1}{3-\varphi}=1,$$
the sequences $V(1,2,0)$ and $V(-1,3,0)$ are complementary in the positive integers. It follows that $R_0=V_0(-1,3,0)$ is the complement of $R_1=V_0(1,2,1)$, by Theorem \ref{th:d0d1}.

\smallskip

\noindent {\rm b)} $w=00$: Theorem \ref{th:d0d1} gives that $R_{00}$ is the union of the two GBS $V_0(1,2,0)$ and $V(3,1,1)$. These two sequences  correspond  to the numbers $N$ with expansions containing $00\cdot 0$, coded $\B$ in \cite{Dekk-phi-FQ}, respectively those  containing $00\cdot 1$, coded $\D$ in \cite{Dekk-phi-FQ}.

\smallskip

\noindent {\rm c)} $w=10$ and  $w=010$: From Theorem \ref{th:d0d1} we obtain that  $R_{10}$ is equal to $V(1,2,-1)$.

\smallskip

\noindent {\rm d)} $w=000$: By Lemma \ref{lem:no} there are no base phi expansions with $d_2d_1d_0d_{-1}(N)=100\cdot 1$. This means that the numbers $N$ from $V(3,1,1)$ in the last part of Theorem \ref{th:d0d1}  do exactly correspond with the numbers $N$ with $d_2d_1d_0d_{-1}(N)=000\cdot 1$. This gives one part of the numbers $N$ where $\beta^+(N)$ has  suffix 000.

The other part comes from the occurrences of $N$ with $d_2d_1d_0d_{-1}(N)=000\cdot 0$. The trick is to observe that the digit blocks
 $1010$ and $000\cdot 0$ always occur in pairs of the expansions of $N-1$ and $N$, for $N=7,\dots 18$. The Propagation Principle (Lemma \ref{lem:prop}, Part {\bf b)})  gives that this coupling will hold for all positive integers $N$.
 From Theorem \ref{th:phi-w0} we know that the digit block $1010$ has occurrence sequence $R_{1010}=V(4,3,-1)$. So the coupling implies that the digit block $000\cdot 0$ has occurrence sequence $V_0(4,3,0)$. Here we should mention that Theorem \ref{th:phi-w0} uses the proposition we are on the way of proving (via the formula $R_{1010}=R_{010}\,B$), however, this only uses part {\rm c)}, which we already proved above.

\smallskip

\noindent {\rm e)} $w=100$:  We already know that expansions with $100\cdot1$ do not occur,  and one checks that an  expansion $\beta(N-2)=\dots 100\cdot0\dots$ always occurs coupled to an expansion $\beta(N)=\dots 00\cdot 1\dots$, for $N=2,\dots,19$. The Propagation Principle (Lemma \ref{lem:prop}, Part {\bf b)})  then implies that this coupling occurs for all $N$. This gives  that $R_{100}=R_{00\cdot 1}-2=V(3,1,-1)$, using the result of part b).
\hfill $\Box$

\bigskip

The sequences $R_{010}$ and $R_{100}$ are examples  of what we call Lucas-Wythoff sequences: their parameters are given respectively by $(L_1,L_0,-1)$ and $(L_2,L_1,-1)$.\\
In general, a \emph{ Lucas-Wythoff sequence} $G$ is a generalized Beatty sequence defined for a natural number $m$ by
$$G = V(L_{m+1},L_m,r),$$
where $r$ is an integer.

\medskip

\begin{theorem}\label{th:phi-w0}
\noindent For any natural number $m\ge 2$ fix a word $w=w_{m-1}\dots w_0$ of $0$'s and $1$'s, containing no $11$. Let $w_0=0$. Then---except if $w=0^m$---the sequence $R_w$ of occurrences of numbers $N$ such that  the  digits $d_{m-1} \dots d_0$ of the base phi expansion of $N$ are equal to $w$, i.e., $d_{m-1}\dots d_0(N)=w$, is a Lucas-Wythoff sequence  of the form
$$R_w= V(L_{m-2},L_{m-3},\gamma_w) \;\;\text{\rm if}\; w_{m-1}=0,\quad R_w=V(L_{m-1},L_{m-2},\gamma_w)  \;\;\text{\rm if}\; w_{m-1}=1,$$
where $\gamma_w$ is a negative integer or 0.\\
In case $w$ consists entirely of $0$'s this sequence of occurrences is given by a disjoint union of two Lucas-Wythoff sequences. We have
\begin{align*}
  R_{0^{2m}} & = V(L_{2m},L_{2m-1},1) \:\sqcup\: V_0(L_{2m-1},L_{2m-2},0),\\
  R_{0^{2m+1}}  & = V_0(L_{2m+1},L_{2m},0) \, \sqcup \, V(L_{2m},L_{2m-1},1).
\end{align*}
\end{theorem}

\noindent{\it Proof:} \:Suppose first that $w$ is a word \emph{not} equal to $0^m$  for some $m\ge 2$.

\noindent The proof is by induction on the length $m$ of $w$. For $m=2$ the statement  of the theorem holds by Proposition \ref{prop:small}, part c).
Next, let $w$ be a word of length $m$ with $w_0=0$.

In the case that $w_{m-1}=1$, $w$ has a unique extension to $0w$, and $R_{0w}=R_w$ is equal to the correct Lucas-Wythoff sequence.

In the case that $w_{m-1}=0$,  the induction hypothesis is  that $R_w$ is a Lucas-Wythoff sequence $R_w=V(L_{m-2},L_{m-3},\gamma_w)$ .

\noindent By Theorem \ref{th:Zeckphi} the numbers $N$ with a $\beta^+(N)$ ending with the digit block $w$ are in one-to-one correspondence with numbers $N'$ with a $Z(N')$ ending with the digit block $w$, and the same property holds for the digit blocks $0w$, respectively $1w$. Note that the correspondence is one-to-one, since the numbers `skipped' in the Zeckendorf expansions all\footnote{We have to follow a different strategy for the words $w=d_{m-1}\dots d_11$ in the next section.} have $d_0=1$. It therefore follows from Proposition 2.6 in \cite{Dekk-Zeck-structure} that $$R_{0w}=R_wA\quad\text{and\:} R_{1w}=R_wB.$$
By Lemma \ref{lem:VA} these have parameters
$$(L_{m-2}+L_{m-3}, L_{m-2}, \gamma_w-L_{m-2})=(L_{m-1},L_{m-2},  \gamma_w-L_{m-2}),$$
respectively
$$(2L_{m-2}+L_{m-1}, L_{m-2}+L_{m-3}, \gamma_w)=(L_m,L_{m-1}, \gamma_w).$$
These are indeed the right expressions for the two words $0w$, respectively $1w$ of length $m+1$.

\medskip

\noindent  Next: The words  $w=0^m$ for some $m\ge 2$.

We claim that for all $m\ge 1$
 \begin{align}
  R_{0^{2m}\cdot 0}     & =  V_0(L_{2m-1},L_{2m-2},0), & \;\;\; R_{0^{2m}\cdot 1} = V(L_{2m},L_{2m-1}, 1)\label{eq:R1} \\
  R_{0^{2m+1}\cdot 0}  & =  V_0(L_{2m+1},L_{2m},0),   & R_{0^{2m+1}\cdot 1}  = V(L_{2m},L_{2m-1}, 1)\label{eq:R2}.
\end{align}
The proof is by induction.

We find in the proof of Proposition  \ref{prop:small}, part b) that $R_{00\cdot 0}=V_0(1,2,0)$ and $R_{00\cdot 1}=V(3,1,1)$.
Since $L_0=2, L_1=1$ and $L_2=3$, this is Equation (\ref{eq:R1}) for $m=1$.

We  find in the proof of Proposition \ref{prop:small}, part d) that $R_{000\cdot 0}=V_0(4,3,0)$ and $R_{000\cdot 1}=V(3,1,1)$. This is Equation (\ref{eq:R2}) for $m=1$.

Next we perform the induction step. Suppose that both Equation (\ref{eq:R1}) and Equation (\ref{eq:R2}) hold.

\smallskip

\noindent {\footnotesize \fbox{(\ref{eq:R1})}} Since $10^{2m+1}\cdot 0$ never occurs by Lemma \ref{lem:no}, we must have
\begin{equation}\label{eq:R0000.0}
 R_{0^{2m+2}\cdot 0}   = R_{0^{2m+1}\cdot 0}  =  V_0(L_{2m+1},L_{2m},0).
 \end{equation}
This is the left part of Equation (\ref{eq:R1}) for $m+1$ instead of $m$.

That $10^{2m+1}\cdot 0$ never occurs also implies that
\begin{equation}\label{eq:R000001.1}
 R_{10^{2m+1}\cdot 1}=R_{10^{2m+1}}=V(L_{2m+1},L_{2m},\gamma_{10^{2m+1}})=V(L_{2m+1},L_{2m}, -L_{2m}+1).
\end{equation}

Here we used the first part of the proof, determining $\gamma_{10^{2m+1}}$ from the observation  that the first occurrence of  $d_{2m+1}\dots d_0(N)=10^{2m+1}$ is at $N=L_{2m+1}+1$, the first element of the Lucas interval $\Lambda_{2m+1}$.

Next we take $V=R_{0^{2m+1}\cdot 1}$, $U=R_{10^{2m+1}\cdot 1}$ and $W=R_{0^{2m+2}\cdot 1}$ in Lemma \ref{lem:VAconv}. According to Equation (\ref{eq:R2}), we take $(p,q,r)=(L_{2m},L_{2m-1},1)$.
The parameters of the sequence $U$ should be $(p+q,p,r-p)=(L_{2m+1},L_{2m},1-L_{2m})$, which conforms with Equation (\ref{eq:R000001.1}).

 The conclusion of Lemma \ref{lem:VAconv} is that $W=R_{0^{2m+2}\cdot 1}$ has parameters
$$(2p+q,p+q,r)=(2L_{2m}+L_{2m-1},L_{2m}+L_{2m-1},1)=(L_{2m+2},L_{2m+1},1).$$
This is the right part of Equation (\ref{eq:R1}) for $m+1$.

\smallskip

\noindent {\footnotesize \fbox{(\ref{eq:R2})}} Since $10^{2m+2}\cdot 1$ never occurs by Lemma \ref{lem:no}, we must have, using the final result of  {\footnotesize\fbox{(\ref{eq:R1})}},
$$ R_{0^{2m+3}\cdot 1}   = R_{0^{2m+2}\cdot 1}  = V( L_{2m+2},L_{2m+1},1).$$
This is the left part of Equation (\ref{eq:R2}) for $m+1$ instead of $m$.

That $10^{2m+2}\cdot 1$ never occurs also implies that
\begin{equation}\label{eq:R10000.0}
 R_{10^{2m+2}\cdot 0}=R_{10^{2m+2}}=V(L_{2m+2},L_{2m+1}, -L_{2m+1}).
\end{equation}

Here we used the first part of the proof, determining $\gamma_{10^{2m+2}}$ from the observation  that the first occurrence of  $d_{2m+3}\dots d_0(N)=10^{2m+2}$ is at $N=L_{2m+2}$, the first element of the Lucas interval $\Lambda_{2m+2}$.

Next we take $V=R_{0^{2m+2}\cdot 0}$, $U=R_{10^{2m+2}\cdot 0}$ and $W=R_{0^{2m+3}\cdot 0}$ in Lemma \ref{lem:VAconv}.
According to Equation (\ref{eq:R0000.0}), we take $(p,q,r)=( L_{2m+1},L_{2m},0)$.
The parameters of the sequence $U$ should be $(p+q,p,r-p)=(L_{2m+2},L_{2m+1},-L_{2m+1})$, which conforms with Equation (\ref{eq:R10000.0}).

 The conclusion of Lemma \ref{lem:VAconv} is that $W=R_{0^{2m+3}\cdot 0}$ has parameters
$$(2p+q,p+q,r)=(2L_{2m+1}+L_{2m},L_{2m+1}+L_{2m},0)=(L_{2m+3},L_{2m+2},0).$$
This is the left part of Equation (\ref{eq:R2}) for $m+1$. \hfill $\Box$

\subsection{Digit blocks $w=d_{m-1}\dots d_11$}\label{sec:w1}

Here there are digit blocks that do not occur at all, like $w=1001$. We denote this as $R_{1001}=\emptyset$.

We order the digit blocks $w$ with $d_0=1$ in a  tree. The first four levels of this tree (taking into account that the node corresponding to $R_{1001}$ has no offspring) are depicted below.

\bigskip

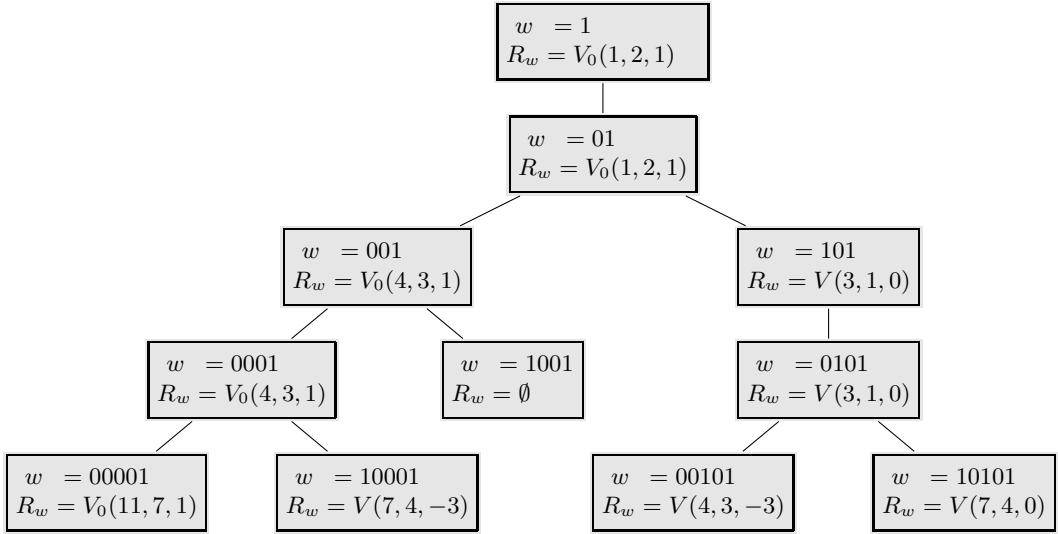
\begin{tikzpicture}
[level distance=15mm,
every node/.style={fill=grey!20,rectangle,inner sep=1pt},
level 1/.style={sibling distance=75mm,nodes={fill=grey!20}},
level 2/.style={sibling distance=60mm,nodes={fill=grey!20}},
level 3/.style={sibling distance=36mm,nodes={fill=grey!20}}]
\node {\footnotesize  $\duoO{=1}{=V_0(1,2,1)\quad}$}            
child {node {\footnotesize $\duoO{=01}{=V_0(1,2,1)}$}            
 child {node {\footnotesize $\duoO{=001}{=V_0(4,3,1)}$}          
  child {node {\footnotesize $\duoO{=0001}{=V_0(4,3,1)}$}       
  child {node {\footnotesize $\duoO{=00001}{=V_0(11,7,1)}$}}       
  child {node {\footnotesize $\duoO{=10001}{=V(7,4,-3)}$}}}       
  child {node {\footnotesize $\duoO{=1001}{=\emptyset}$}           
    }     
}
 child {node {\footnotesize $\duoO{=101}{=V(3,1,0)}$}         
  child {node {\footnotesize $\duoO{=0101}{=V(3,1,0)}$}     
  child {node {\footnotesize $\duoO{=00101}{=V(4,3,-3)}$}}       
  child {node {\footnotesize $\duoO{=10101}{=V(7,4,0)}$}}}       
} };
\end{tikzpicture}

\bigskip

Here $R_{01}=R_1=V_0(1,2,1)$ has been given in Theorem \ref{th:d0d1}. The correctness of the other occurrence sequences follows  from Theorem \ref{th:phi-w1}.

We next determine an infinite family of excluded blocks.

\medskip

\begin{lemma}\label{lem:1001} Let $m\ge 2$ be an integer. There are no expansions $\beta^+(N)$ with end block $10^{2m}1$.
\end{lemma}

\noindent{\it Proof:} Consider any $N$ such that $\beta^+(N)$ has end block $10^{2m}1$. Such an $N$, of course, has $d_{-1}(N)=0$, so we see that
$\beta(N-1)=\dots 10^{2m+1}\cdot 0\dots$. According to Lemma \ref{lem:no} this is not possible. \hfill $\Box$

\medskip
Next, we establish a connection with the previous section.

\begin{lemma} \label{lem:1to0} Let $m\ge 2$ be an integer. The  block $w=d_{m-1}\dots d_11\cdot 0$ is end block of  $\beta^+(N)$ if and only if
the block $\breve{w}:=d_{m-1}\dots d_10\cdot 0$ occurs in $\beta(N-1)$. 
\end{lemma}

\noindent{\it Proof:} This follows quickly from the Propagation Principle Lemma \ref{lem:prop} applied to the couple of blocks $00\cdot 0$ and $01\cdot 0$. \hfill $\Box$

\begin{theorem}\label{th:phi-w1}
\noindent For any natural number $m\ge 2$ fix a word $w=w_{m-1}\dots w_0$ of $0$'s and $1$'s, containing no $11$. Let $w_0=1$.
With exception of the words $w$ with suffix $0^m1$ and $10^m1$, for $m=2,3,\dots$, the sequence $R_w$ of occurrences of numbers $N$ such that  the  digits $d_{m-1} \dots d_0$ of the base phi expansion of $N$ are equal to $w$, i.e., $d_{m-1}\dots d_0(N)=w$, is a Lucas-Wythoff sequence  of the form
$$R_w= V(L_{m-2},L_{m-3}, \gamma_w) \;\;\text{\rm if}\; w_{m-1}=0,\quad R_w=V(L_{m-1},L_{m-2}, \gamma_w)  \;\;\text{\rm if}\; w_{m-1}=1,$$
where $\gamma_w$ is a negative integer or 0.

In case $w=0^{2m}1$ we have $R_w=V_0(L_{2m+1},L_{2m},1)$, and this is also the sequence of occurrences of   $w=0^{2m+1}1$.

In case $w=10^{2m}1$ the word $w$ does not occur at all as digit end block.

In case $w=10^{2m+1}1$ we have $R_w=V(L_{2m+2},L_{2m+1},-L_{2m+1}+1)$.
\end{theorem}

\noindent{\it Proof:} 
It follows from Lemma \ref{lem:1to0} that $R_w=R_{\breve{w}}+1$, if $R_w\ne \emptyset$. So the first part of Theorem \ref{th:phi-w0} yields the statement of the theorem for all $w$ not equal to $0^m1$ or $10^m1$.

In case $w=0^{2m}1\cdot 0$, we have  $\breve{w}=0^{2m+1}\cdot 0$, and the result follows from the left part of Equation (\ref{eq:R2}).

In case $w=10^{2m}1$ the word $w$ does not occur as digit end block, according to Lemma \ref{lem:no}.

In case $w=10^{2m+1}1\cdot 0$ we have   $\breve{w}=10^{2m+2}\cdot 0$,  and now Equation (\ref{eq:R10000.0})  gives that
           $R_w=R_{\breve{w}}+1=V(L_{2m+2},L_{2m+1},-L_{2m+1}+1)$.\hfill $\Box$

\section{The negative powers of the golden mean}\label{sec:neg}

Here we discuss what we can say about the words $\beta^-(N)$.
These do  have an even more intricate structure than the $\beta^+(N)$.

\subsection{The words $\beta^-(N)$}\label{sec:trident}

Here we  look at complete $\beta^-(N)$'s. Although at first sight these seem to appear in a random order, there is an order dictated not by a coin toss, but by another dynamical system: the rotation over an angle $\varphi$. Moreover, they appear in singletons, or as triples. This can be proved with the \{$\A\B\C$, $\D$\}--structure found in the paper \cite{Dekk-phi-FQ}.

For a more extensive analysis, partition the natural numbers larger than 1 into intervals
$$\Xi_n:=\Lambda_{2n-1}\cup\Lambda_{2n}=[L_{2n-1}+1,L_{2n+1}].$$
The relevance of the $\Xi_n, n=1,2,\dots$ is that these are exactly the intervals where $\beta^-(N)$ has length $2n$.
The  $\Xi_n$ intervals have length $$L_{2n+1}-L_{2n-1}=L_{2n+1}-L_{2n}+L_{2n}-L_{2n-1}=L_{2n-1}+L_{2n-2}=L_{2n}.$$
Call three consecutive numbers $N,N+1,N+2$ a \emph{trident}, if $\beta^-(N)=\beta^-(N+1)=\beta^-(N+2)$. For example: 2,3,4 and 6,7,8  are tridents.
We shall always take the middle number $N\!+\!1$ as the representing number of a trident interval $[N,N\!+\!1,N\!+2]$. We call this number $\Pi$-\emph{essential}.
By definition the other $\Pi$-essential numbers are the singletons.

\begin{lemma} {\bf [Trident splitting]} \label{lem:trident}
In $\Lambda_{2n-1}\cup\Lambda_{2n}$ the last number  of $\Lambda_{2n-1}$ and the first two numbers in $\Lambda_{2n}$ are in
 the same trident.
\end{lemma}

\noindent {\it Proof:} This is true for $n=1$ and $n=2$: $\Lambda_1\cup\Lambda_2=\{2\}\cup[3,4]$ is a trident,
and $\Lambda_3\cup\Lambda_4=[5,6]\cup[7,8,\dots,11]$ contains the trident $[6,7,8]$.
The property then follows by induction, using Theorem \ref{th:rec2}. \hfill $\Box$

\medskip

The following lemma helps to count singletons and  tridents.

\begin{lemma} The following relation between Lucas numbers and Fibonacci numbers holds: $F_n+3F_{n+1}=L_{n+2}$ for $n=0,1,2,\dots$.
\end{lemma}

For a proof, note that $F_0+3F_1=3=L_2$, and $F_2+3F_3=1+6=L_4$, and then add these two equations, etc.

The lemma describes the fact that the $\Xi_n$ intervals contain $F_{2n-2}$ singletons, and $F_{2n-1}$ tridents, making a total number of $L_{2n}$.
The collection of different $\beta^-(N)$-blocks of length $2n$  has thus cardinality $F_{2n-2}+F_{2n-1}=F_{2n}$. This implies that we have proved the following theorem.

\begin{theorem}\label{th:all-beta-min} All Zeckendorf words of even length ending in 1 appear as $\beta^-(N)$-blocks.
\end{theorem}


Here we mean by a Zeckendorf word (or golden mean word) all words in which 11 does not occur. We denote by $\mathcal{Z}_{m}$ the set of Zeckendorf words of length $m$, for $m=1,2,\dots$.  It is easily proved that the cardinality of $\mathcal{Z}_{m}$ equals $F_{m+2}$.
So the cardinality of the set of  words from $\mathcal{Z}_{2n}$ ending in 1 is equal to $F_{2n}$, implying the result of Theorem \ref{th:all-beta-min}.


\bigskip

 Since all   $\beta^-(N)$ have suffix 01, the essential information of these words is contained in
 $$    \gamma^-(N):=  \beta^-(N)1^{-1}0^{-1}.$$
 The words   $\gamma^-(N)$ are Zeckendorf words, corresponding one-to-one to the natural numbers  $Z^{-1}(\gamma^-(N))$.
 Obviously, the $\gamma^-(N)$ have the same ordering as the $\beta^-(N)$.
 According to Theorem \ref{th:all-beta-min} we then (after identifying tridents with their middle number) obtain a permutation of length $F_{2n}$  of the $\Pi$-essential elements of $\Xi_n$ by  coding these numbers by  $\C(N):=Z^{-1}(\gamma^-(N))$. We  denote this permutation by $\Pi^\beta_{2n}$.

 The following Zeckendorf words and codes will be important in the sequel.

 \begin{lemma} \label{lem:code} For all natural numbers $n$ we have
 \begin{equation}\label{eq:bordergammas}
    \gamma^-(L_{2n})=0^{2n-2}, \;  \gamma^-(L_{2n+1})=[01]^{n-1}, \; \gamma^-(L_{2n+1}+1)=[10]^n,\; \gamma^-(L_{2n+2}-1)=0^{2n}\!.
 \end{equation}
 \begin{equation}\label{eq:bordercodes}
   \C(L_{2n})=0, \quad  \C(L_{2n+1})=F_{2n-1}-1, \quad  \C(L_{2n+1}+1)=F_{2n+2}-1, \quad  \C(L_{2n+2}-1)=0.
 \end{equation}
 \end{lemma}

\noindent {\it Proof:}
The correctness of Equation (\ref{eq:bordergammas}) follows from Equations (\ref{eq:Lm}) and (\ref{eq:Lmplus1}).
So  $\gamma^-(L_{2n})$ is the first word in $\mathcal{Z}_{2n-2}$,  $\gamma^-(L_{2n+1})$ is 0  followed by the last word in  $\mathcal{Z}_{2n-3}$,
$\gamma^-(L_{2n+1}+1)$ is the last word in $\mathcal{Z}_{2n}$, and
$\gamma^-(L_{2n+2}-1)$ is the first word in  $\mathcal{Z}_{2n-2}$.
Since $\mathcal{Z}_{m}$ has cardinality $F_{m+2}$, Equation (\ref{eq:bordercodes}) follows.  \hfill $\Box$

\medskip

 We have to determine the codings of all natural numbers $N$. For this it is useful to translate Theorem \ref{th:rec2} to the $\gamma^-\!$-blocks.

 \begin{theorem}{\bf [Recursive structure theorem: $\gamma^-\!$-version]}\label{th:recg}\\
\noindent{\,\bf (i):  Odd\;} For all $n\ge 1$ one has $\Lambda_{2n+1}=\Lambda^{(a)}_{2n-1}\cup\Lambda^{(b)}_{2n-2}\cup\Lambda^{(c)}_{2n-1}, $
where $\Lambda^{(a)}_{2n-1}=\Lambda_{2n-1}+L_{2n}$,\; $\Lambda^{(b)}_{2n-2}=\Lambda_{2n-2}+L_{2n+1}$, and $\Lambda^{(c)}_{2n-1}=\Lambda_{2n-1}+L_{2n+1}$.\\
We have\\[-.8cm]
\begin{subequations} \label{eq:shift-odd}
\begin{align}
  \gamma^-(N)= & \;\gamma^-(N-L_{2n})\,10& for\; N\in \Lambda^{(a)}_{2n-1}, \label{eq:15a}\\
  \gamma^-(N)= & \;\gamma^-(N-L_{2n+1})\,0010&      for\; N\in  \Lambda^{(b)}_{2n-2},\label{eq:15b}\\
  \gamma^-(N)= & \;\gamma^-(N-L_{2n+1})\,00 &         for\; N\in \Lambda^{(c)}_{2n-1}\label{eq:15c}.
\end{align}\\[-.6cm]
\end{subequations}
\noindent{\,\bf (ii): Even\;} For all $n\ge 1$ one has $\Lambda_{2n+2}=\Lambda^{(a)}_{2n}\cup\Lambda^{(b)}_{2n-1}\cup\Lambda^{(c)}_{2n}, $
 where $\Lambda^{(a)}_{2n}=\Lambda_{2n}+L_{2n+1}$,\; $\Lambda^{(b)}_{2n-1}=\Lambda_{2n-1}+L_{2n+2}$, and $\Lambda^{(c)}_{2n}=\Lambda_{2n}+L_{2n+2}$.\\
 We have\\[-.8cm]
\begin{subequations} \label{eq:shift-even}
\begin{align}
  \gamma^-(N)= & \;\gamma^-(N-L_{2n+1})\,00 & for\; N\in \Lambda^{(a)}_{2n},\phantom{x} \label{eq:16a} \\
  \gamma^-(N)= & \; \gamma^-(N-L_{2n+2})\,01 &                     for\; N\in \Lambda^{(b)}_{2n-1},  \label{eq:16b}\\
  \gamma^-(N)= & \;\gamma^-(N-L_{2n+1})\,01 &                       for\; N\in \Lambda^{(c)}_{2n}.\phantom{x.}  \label{eq:16c}
\end{align}
\end{subequations}
\end{theorem}
                                                                                                                  \!
 We give the situation for $n=2$, where $\Xi_2=\Lambda_3 \cup \Lambda_4=[5,6,\dots,11]$.

\bigskip

 \begin{tabular}{|c|c|c|c|c|}
   \hline
  \; $N^{\phantom{|}}$ & $\Lambda$-int. & $\cdot\beta^-(N)$    & $\cdot\gamma^-(N)$    & $\C(N)$  \\[.0cm]
   \hline
   5\;    & $\Lambda_3$ & $\cdot1001$        & \; $\cdot10$    & \; 2   \\
   6\;    & $\Lambda_3$ & $\cdot0001$        & \; $\cdot00$    & \; \grijs{0}    \\
   \hline
   7\;    & $\Lambda_4$ & $\cdot0001$         & \; $\cdot00$    & \;      0     \\
   8\;    & $\Lambda_4$ & $\cdot0001$         & \; $\cdot00$    & \; \grijs{0}     \\
   9\;    & $\Lambda_4$ & $\cdot0101$         & \; $\cdot01$    & \; \grijs{1}     \\
   10\,   & $\Lambda_4$ & $\cdot0101$         & \; $\cdot01$    & \; 1     \\
   11\,   & $\Lambda_4$ & $\cdot0101$        & \; $\cdot01$    & \; \grijs{1}     \\
   \hline
 \end{tabular}

 \medskip

\noindent We see that  $\Pi^\beta_{4}=\big( 2\, 0\, 1\big)$.

  Here is the situation for $n=3$, where $\Xi_3=\Lambda_5 \cup \Lambda_6=[12,13,\dots,29]$.

 \bigskip

 \begin{tabular}{|c|c|c|c|c|}
   \hline
  \; $N^{\phantom{|}}$ &{\small  $\Lambda$-int.}&  $\cdot\beta^-(N)$ & $\cdot\gamma^-(N)$ & $\C(N)$  \\[.0cm]
   \hline
   12\; &  $\Lambda_5$ &  $\cdot101001$  &  $\cdot1010$ &  7  \\
   13\; &  $\Lambda_5$ &  $\cdot001001$  &  $\cdot0010$ &  \grijs{2}  \\
   14\; &  $\Lambda_5$ &  $\cdot001001$  &  $\cdot0010$ &  2  \\
   15\; &  $\Lambda_5$ &  $\cdot001001$  &  $\cdot0010$ &  \grijs{2}  \\
   16\; &  $\Lambda_5$ &  $\cdot100001$  &  $\cdot1000$ &  5  \\
   17\; &  $\Lambda_5$ &  $\cdot000001$  &  $\cdot0000$ & \grijs{0}  \\
   \hline
   18\; &  $\Lambda_6$ &  $\cdot000001$  &  $\cdot0000$ &  0  \\
   19\; &  $\Lambda_6$ &  $\cdot000001$  &  $\cdot0000$ &  \grijs{0}  \\
   20\; &  $\Lambda_6$ &  $\cdot010001$  &  $\cdot0100$ &  \grijs{3}  \\
   \hline
 \end{tabular}\quad
 \begin{tabular}{|c|c|c|c|c|}
   \hline
  \; $N^{\phantom{|}}$ &{\small  $\Lambda$-int.}&  $\cdot\beta^-(N)$ & $\cdot\gamma^-(N)$ & $\C(N)$  \\[.0cm]
   \hline
   21\; &  $\Lambda_6$ &  $\cdot010001$  &  $\cdot0100$ &  3  \\
   22\; &  $\Lambda_6$ &  $\cdot010001$  &  $\cdot0100$ &  \grijs{3}  \\
   23\; &  $\Lambda_6$ &  $\cdot100101$  &  $\cdot1001$ &  6  \\
   24\; &  $\Lambda_6$ &  $\cdot000101$  &  $\cdot0001$ &  \grijs{1}  \\
   25\; &  $\Lambda_6$ &  $\cdot000101$  &  $\cdot0001$ &  1  \\
   26\; &  $\Lambda_6$ &  $\cdot000101$  &  $\cdot0001$ &  \grijs{1}  \\
   27\; &  $\Lambda_6$ &  $\cdot010101$  &  $\cdot0101$ &  \grijs{4}  \\
   28\; &  $\Lambda_6$ &  $\cdot010101$  &  $\cdot0101$ &  4  \\
   29\; &  $\Lambda_6$ &  $\cdot010101$  &  $\cdot0101$ &  \grijs{4}  \\
   \hline
 \end{tabular}\quad

 \medskip

\noindent We see that $\Pi^\beta_{6}=\big(7\,2\,5\,0\,3\,6\,1\,4\big)$.

\medskip

What are these permutations?

\begin{theorem}\label{th:permut} For all natural numbers $n$  consider the $F_{2n}$ Zeckendorf words of length $2n$  occurring as $\beta^-(N)$ in the $\beta$-expansions of the numbers in $\Xi_n$. Then these occur in an order given by a permutation $\Pi^\beta_{2n}$ which is the  orbit of the element $F_{2n}-1$ under the addition by the element $F_{2n-2}$ on the cyclic group $\mathbb{Z}/F_{2n}\mathbb{Z}$.
\end{theorem}

\noindent {\it Proof:} We  have to show for all $n$ that
\begin{equation}\label{eq:IH}
  \Pi^\beta_{2n}(1)=F_{2n}-1,\quad \Pi^\beta_{2n}(j+1)=\Pi^\beta_{2n}(j)+F_{2n-2} \!\mod F_{2n}, \; {\rm for}\;j=1,\dots,F_{2n}-1.
\end{equation}
It is easily checked that the cases $n=2$ and $n=3$ given above conform with this. For $n=3$ one has: $F_6=8$, $F_4=3$, and
 $\Pi^\beta_{6}(1)=7, \,\Pi^\beta_{6}(j+1) =\Pi^\beta_{6}(j)+3 \mod 8$ for $j=1,\dots 7$.

 \medskip

The first claim in Equation (\ref{eq:IH}) follows from Lemma \ref{lem:code} for all $n$: since the interval $\Xi_n=[L_{2n-1}+1,L_{2n+1}]$, we have $\Pi^\beta_{2n}(1)=F_{2n}-1$ according to  Equation (\ref{eq:bordercodes}).


\medskip

The proof proceeds by induction, based on Theorem \ref{th:recg}, the $\gamma^-\!$-version of the Recursive Structure Theorem.

For the second part of Equation (\ref{eq:IH}) with $n$ replaced by $n+1$, we have to split the permutation $\Pi_{2n+2}^\beta$ into six pieces,
and then we have to glue the expressions together to obtain
the full permutation on  the set $\Xi_{n+1} = \Lambda_{2n+1} \cup \Lambda_{2n+2} = [L_{2n+1}+1, L_{2n+2}-1] \cup [L_{2n+2},L_{2n+3}]$.
According to the Recursive Structure Theorem
\begin{equation}\label{eq:six}
\Xi_{n+1} =  \Lambda^{(a)}_{2n-1}\cup\Lambda^{(b)}_{2n-2}\cup\Lambda^{(c)}_{2n-1}\cup\Lambda^{(a)}_{2n}\cup\Lambda^{(b)}_{2n-1}\cup\Lambda^{(c)}_{2n}.
\end{equation}
We start with the first interval, $\Lambda^{(a)}_{2n-1}$. From Theorem \ref{th:recg} we have that for $N\in \Lambda^{(a)}_{2n-1}$,
\begin{equation}\label{eq:grec1}
\gamma^-(N)= \gamma^-(N-L_{2n})\,10.
\end{equation}
 What does this imply for the codes?

 Let $Z(\C(N-L_{2n}))=\gamma^-(N-L_{2n})=d_{2n-3}\dots d_0$, so $\C(N-L_{2n})= \sum_{i=0}^{2n-3} d_i \ddot{F}_i$. Then Equation (\ref{eq:grec1}) leads to
 $$   \C(N)=\sum_{i=0}^{2n-3} d_i \ddot{F}_{i+2} + 1\cdot \ddot{F}_1 + 0\cdot \ddot{F}_0=\sum_{i=0}^{2n-3} d_i \ddot{F}_{i+2}+2.$$
This implies, in particular,   that the differences between the codes of two consecutive $\Pi$-essential numbers within  the interval $\Lambda_{2n-1}$  have increased from $F_{2n-2}\!\mod F_{2n}$ to $F_{2n} \!\mod F_{2n+2}$ for the corresponding numbers in the interval $\Lambda^{(a)}_{2n-1}$.

 We pass to the second  interval, $\Lambda^{(b)}_{2n-2}$. From Theorem \ref{th:recg} we have that for $N$ from $\Lambda^{(b)}_{2n-2}$,
\begin{equation}\label{eq:grec2}
\gamma^-(N)= \gamma^-(N-L_{2n+1})\,0010.
\end{equation}
 What does this imply for the codes?

 Let $Z(\C(N-L_{2n+1}))=\gamma^-(N-L_{2n+1})=d_{2n-4}\dots d_0$, so  $\C(N-L_{2n+1})= \sum_{i=0}^{2n-4} d_i \ddot{F}_i$. Then Equation (\ref{eq:grec2}) leads to
 $$\C(N)=\sum_{i=0}^{2n-4} d_i \ddot{F}_{i+4} + 0\cdot\ddot{F}_3 + 0\cdot \ddot{F}_2 +1\cdot \ddot{F}_1 + 0\cdot \ddot{F}_0=\sum_{i=0}^{2n-4} d_i \ddot{F}_{i+4}+2.$$
This implies  that the differences between the codes of two consecutive numbers within  the interval $\Lambda_{2n-2}$
 have increased from $F_{2n-4}\!\mod F_{2n-2}$ to $F_{2n} \!\mod F_{2n+2}$ for the corresponding numbers in the interval $\Lambda^{(b)}_{2n-2}$.

 Similar computations give that for the next 4 intervals $\Lambda^{(c)}_{2n-1}, \Lambda^{(a)}_{2n},\Lambda^{(b)}_{2n-1}$, and $\Lambda^{(c)}_{2n}$ there always is an addition of $F_{2n} \!\mod F_{2n+2}$.

 \medskip

 The remaining task is to check that the same holds on the five boundaries between the translated $\Lambda$-intervals.
 We number these boundaries with the roman numerals I,  II, III, IV, V.

 \medskip

\noindent \fbox{III \& V:}\;For the  third and the fifth boundary between respectively  the intervals $\Lambda^{(c)}_{2n-1}$ and $\Lambda^{(a)}_{2n}$ and the intervals $\Lambda^{(b)}_{2n-1}$ and $\Lambda^{(c)}_{2n}$ this follows from the Trident Splitting Lemma, Lemma \ref{lem:trident}.
 The  reason is that if $[N,N+1,N+2]$ is the trident which is splitted, then the difference between $\C(N-1)$ and $\C(N)$ is equal to $F_{2n} \!\mod F_{2n+2}$, as these two numbers are both from the first translated $\Lambda$-interval, and not from the same trident. But then  the difference between the codes of the last $\Pi$-essential number $N-1$ in the first translated $\Lambda$-interval, and the first $\Pi$-essential number $N+1$ in the second translated $\Lambda$-interval is also equal to $F_{2n} \!\mod F_{2n+2}$.

 \medskip

\noindent \fbox{ I:}\;   The last number in the first interval $\Lambda^{(a)}_{2n-1}$ is $2L_{2n}-1$ with associated $\gamma^-\!$-block
 $$\gamma^-(2L_{2n}-1)= \gamma^-(2L_{2n}-1-L_{2n})\,10=\gamma^-(L_{2n}-1)\,10=0^{2n-1}\,10.$$
 Here we used Equation (\ref{eq:shift-odd}a) in the first, and  Equation (\ref{eq:Lmplus1}) in the last step.   It follows directly that $\C(2L_{2n}-1)=2$.

 The first number in the second interval $\Lambda^{(b)}_{2n-2}$ is $2L_{2n}$. From Equation (\ref{eq:Lm}) we have  $\beta(2L_{2n})\doteq 20^{2n}\cdot 0^{2n-1}2\doteq 20^{2n}\cdot 0^{2n-1}1001$, so $\gamma^-(2L_{2n})=0^{2n-1}10$, giving  $\C(2L_{2n})=2$. It is clear that also the second number $2L_{2n}+1$ in $\Lambda^{(b)}_{2n-2}$ has code  $\C(2L_{2n}+1)=2$. As in the previous case, this implies that the difference between the codes of the last $\Pi$-essential number in the first translated $\Lambda$-interval, and the first $\Pi$-essential number in the second translated $\Lambda$-interval is equal to $F_{2n} \!\mod F_{2n+2}$.

 \medskip

\noindent \fbox{ II:}\; The last number in the second interval $\Lambda^{(b)}_{2n-2}$  is the number $L_{2n-1}+L_{2n+1}$.  According to Equation (\ref{eq:shift-odd}b) the associated $\gamma^-\!$-block is
 $$\gamma^-(L_{2n-1}+L_{2n+1})= \gamma^-(L_{2n-1}+L_{2n+1}-L_{2n+1})\,0010=\gamma^-(L_{2n-1})\,0010=[01]^{n-2}\,0010.$$
But we know from Lemma \ref{lem:code} that \; $\gamma^-(L_{2n-1})\,0101=[01]^n=\gamma^-(L_{2n+3}).$

 By Lemma \ref{lem:code}  we have that $\C(L_{2n+3})=F_{2n+1}-1$. To obtain the code of $N=L_{2n-1}+L_{2n+1}$, we have to subtract the number $F_3+F_1=3$ with Zeckendorf expansion $0101$, and add the number $F_2=2$  with Zeckendorf expansion $0010$. This gives the code
 $$\C(L_{2n-1}+L_{2n+1})=F_{2n+1}-1-3+1=F_{2n+1}-3.$$
The first number in the third interval $\Lambda^{(c)}_{2n-1}$ is the number $L_{2n-1}+L_{2n+1}+1$. According to   according to  Equation (\ref{eq:shift-odd}c) the associated $\gamma^-\!$-block is
 $$\gamma^-(L_{2n-1}+L_{2n+1}+1)= \gamma^-(L_{2n-1}+L_{2n+1}+1-L_{2n+1})\,00=\gamma^-(L_{2n-1}+1)\,00.$$
But we know from Lemma \ref{lem:code} that \; $\gamma^-(L_{2n-1}+1)10=[10]^n=\gamma^-(L_{2n+1}+1).$

By Lemma \ref{lem:code} we have that $\C(L_{2n+1}+1)=F_{2n+2}-1$.
To obtain the code of $N=L_{2n-1}+L_{2n+1}+1$, we have to subtract the number $F_2=2$ with Zeckendorf expansion $10$, from this code. This gives the code
 $$\C(L_{2n-1}+L_{2n+1}+1)=F_{2n+2}-1-2=F_{2n+2}-3.$$
The  conclusion is that $L_{2n-1}+L_{2n+1}$ and $N=L_{2n-1}+L_{2n+1}+1$ are $\Pi$-essential, with difference in codes  $F_{2n+2}-3-(F_{2n+1}-3)=F_{2n}.$

\medskip

\noindent \fbox{ IV:}\;  The last number in the fourth interval $\Lambda^{(c)}_{2n}$  is the number $L_{2n+1}+L_{2n+1}=2L_{2n+1}$. According to Equation (\ref{eq:shift-even}a) the associated $\gamma^-\!$-block is
 $$\gamma^-(2L_{2n+1})= \gamma^-(2L_{2n+1}-L_{2n+1})\,00=\gamma^-(L_{2n+1})\,00=[01]^{n-1}\,00.$$
But we know from Lemma \ref{lem:code} that \; $\gamma^-(L_{2n+1})\,01=[01]^n=\gamma^-(L_{2n+3}).$

By Lemma \ref{lem:code}  we have that $\C(L_{2n+3})=F_{2n+1}-1$. To obtain the code of $N=2L_{2n+1}$, we have to subtract the number $F_1=1$ with Zeckendorf expansion $01$. This gives the code
 $$\C(2L_{2n+1})=F_{2n+1})-1-1=F_{2n+1}-2.$$

 The first number in the fifth interval $\Lambda^{(b)}_{2n-1}$ is the number $L_{2n-1}+1+L_{2n+2}$. According to Equation (\ref{eq:shift-even}b) the  associated $\gamma^-\!$-block is
 $$\gamma^-(L_{2n-1}+1+L_{2n+2})= \gamma^-(L_{2n-1}+1+L_{2n+2}-L_{2n+2})\,01=\gamma^-(L_{2n-1}+1)\,01.$$
 But we know from Lemma \ref{lem:code} that \; $\gamma^-(L_{2n-1}+1)10=[10]^n=\gamma^-(L_{2n+1}+1).$

 By Lemma \ref{lem:code} we have that $\C(L_{2n+1}+1)=F_{2n+2}-1$.
 To obtain the code of $N=L_{2n-1}+1+L_{2n+2}$, we have to subtract the number $F_2=2$ with Zeckendorf expansion $10$, and add the number $F_1=1$ with Zeckendorf expansion 01 to this code. This gives the code
 $$\C(L_{2n-1}+L_{2n+1}+1)=F_{2n+2}-1-2+1=F_{2n+2}-2.$$
The  conclusion is that $2L_{2n+1}$ and $L_{2n-1}+1+L_{2n+2}$ are $\Pi$-essential, with difference in codes  $F_{2n+2}-2-(F_{2n+1}-2)=F_{2n}.$ \hfill $\Box$

\bigskip

We now explain the connection with a rotation on a circle mentioned at the beginning of this section. Note that with this point of view all the cyclic groups of Theorem \ref{th:permut} are represented by a single object: the rotation on the circle.

\begin{theorem}\label{th:rot} For all natural numbers $n$ the permutations $\Pi^\beta_{2n}$ are given by the order in which  the first $F_{2n}$ iterates of the rotation $z\rightarrow \exp(2\pi i (z-\varphi))$ occur on the circle.
\end{theorem}

We sketch a proof of this result based on the paper  \cite{Ravenstein-1988}. In the literature one will not find the rotation $z\rightarrow \exp(2\pi i (z-\varphi))$, but several papers treat  the rotation $z\rightarrow \exp(2\pi i (z+\tau))$, where $\tau$ is the algebraic conjugate of $\varphi$. Note that this rotation has exactly the same orbits as $z\rightarrow \exp(2\pi i (z+\varphi))$, and replacing $\varphi$ by $-\varphi$ amounts to reversing the permutation. In the literature the origin is usually added to the orbit. For instance in \cite{Ravenstein-1988}, the $N$ ordered iterates are given by the permutation $\big(u_1\,u_2\,\dots \, u_N\big)$, which for \emph{all} $N$ gives a permutation starting trivially with $u_1=0$.

Lemma 2.1 in \cite{Ravenstein-1988} states that for $j=1,\dots,N$ one has $u_j=(j-1)u_2 \mod N$.\\
Next, Theorem 3.3 in \cite{Ravenstein-1988} states that $u_2=u_2(N)= F_{2n-1}$  in the case that $N=F_{2n}$, $n\ge 1$.

We illustrate this for the case $n=3$.\\   
We  have $N=F_6=8$, and $0< \{5\tau\}<\{2 \tau\}<\{ 7\tau\}<\{4\tau\}<\{\tau\}<\{6 \tau\}<\{ 3\tau\},$
so $ \big(u_1\,u_2\,\dots \, u_N\big) = \big(0\,5\,2\,7\,4\,1\,6\,3\big)$. 
As $\{8\tau\}$ is the largest number in the rotation orbit of the first 9 iterations,
$\big(u_{N+1}\,u_N\,\dots \, u_2\big)=\big(8\,3\,6\,1\,4\,7\,2\,5\big).$
After subtraction of 1 in all entries, one obtains the permutation $\Pi^\beta_{6}$.

\subsection{Digit blocks  $w=d_{-1}\dots d_{-m}(N)$ as prefix of $\beta^-(N)$}

  For any digit block $w$  we will try to determine the  sequence $R_w$ of those numbers $N$ with  $w$ as prefix of $\beta^-(N)$.
  The tridents introduced in the previous section give occurrence sequences $R_w$ which are unions of three consecutive generalized Beatty sequences.
  We will write for short
  $$V(p,q,[r,r+1,r+2]):=V(p,q,r)\sqcup V(p,q,r+1)\sqcup V(p,q,r+2).$$
  As before, we order the  $w$  in a Fibonacci tree. Here we write $R_{\cdot w}$ for the occurence sequences of words $w$ occurring as a prefix of the words $\beta^-(N)$, to emphasize the positions of these words in the expansion $\beta(N)$. The first four levels of this tree are depicted below.

\bigskip

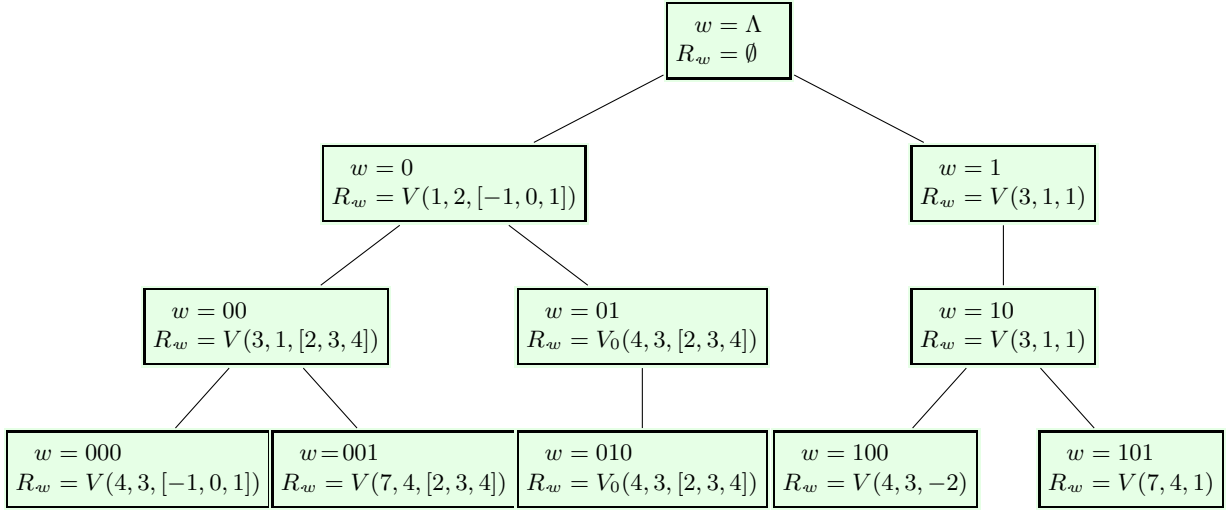
\hspace*{-0.8cm}\begin{tikzpicture}
[level distance=19mm,
every node/.style={fill=green!10,rectangle,inner sep=1pt},
level 1/.style={sibling distance=73mm,nodes={fill=green!10}},
level 2/.style={sibling distance=50mm,nodes={fill=green!10}},
level 3/.style={sibling distance=34mm,nodes={fill=green!10}}]
\node {\footnotesize  $\duoZ{=\Lambda}{=\emptyset\quad}$}            
child {node {\footnotesize $\duoZ{=0}{=V(1,2,[-1,0,1])}$}            
 child {node {\footnotesize $\duoZ{=00}{=V(3,1,[2,3,4])}$}          
  child {node {\footnotesize $\duoZ{=000}{=V(4,3,[-1,0,1])}$}}       
  child {node {\footnotesize $\duoZ{\!=\! 001}{=V(7,4,[2,3,4])}$}}    
}
 child {node {\footnotesize $\duoZ{=01}{=V_0(4,3,[2,3,4])}$}         
  child {node {\footnotesize $\duoZ{=010}{=V_0(4,3,[2,3,4])}$}}     
}
}
child {node {\footnotesize $\duoZ{=1}{=V(3,1,1)}$}            
 child {node {\footnotesize $\duoZ{=10}{=V(3,1,1)}$}          
  child {node {\footnotesize $\duoZ{=100}{=V(4,3,-2)}$}}     
   child {node {\footnotesize $\duoZ{=101}{=V(7,4,1)}$}}  
}
};
\end{tikzpicture}

\bigskip

We start with the words $w$ on this tree.

\medskip

\begin{proposition}\label{prop:smallminus} Let $\beta(N)=\beta^+(N)\cdot\beta^-(N)$ be the base phi expansion of the   number $N$.\\
  Let $w$ be a word of length $m$. Then the sequence of occurrences $R_w$ of numbers $N$ such that  the  first $m$ digits of $\beta^-(N)$  are equal to $w$, i.e., $d_{-1}\dots d_{-m}(N)=w$, is given for  the words $w$ of length at most 3, by\\
\mbox{\rm a)} $R_{\cdot 0} = V(2,1,-1)\, \sqcup \,V(2,1,0) \, \sqcup \,V(2,1,1)$,\\
\mbox{\rm b)} $R_{\cdot 1}= R_{\cdot 10} = V(3,1,1)$,\\
\mbox{\rm c)} $R_{\cdot 00} = V(3,1,2)\, \sqcup \,V(3,1,3) \, \sqcup \,V(3,1,4)$,\\
\mbox{\rm d)} $R_{\cdot 01} =R_{\cdot 010}= V_0(4,3,2)\, \sqcup\, V_0(4,3,3)\, \sqcup \,V_0(4,3,4)$,\\
\mbox{\rm e)} $R_{\cdot 000} = V(4,3,-1)\, \sqcup\, V(4,3,0)\, \sqcup \,V(4,3,1)$,\\
\mbox{\rm f)} $R_{\cdot 001} = V(7,4,2)\, \sqcup \, V(7,4,3)\, \sqcup \, V(7,4,4)$,\\
\mbox{\rm g)} $R_{\cdot 100} = V(4,3,-2)$,\\
\mbox{\rm h)} $R_{\cdot 101} = V(7,4,1)$.
\end{proposition}

\medskip

\noindent {\it Proof:}

 \noindent {\rm a)} $w=\cdot 0$: In Section 5 of the paper \cite{Dekk-phi-FQ} the tridents are coded by triples $(\A, \B, \C)$. It follows from Theorem 5.1 of \cite{Dekk-phi-FQ} that the first elements (coded $\A$) of the tridents are all member of $V(2,1,-1)$. This implies the statement in {\rm a)}.

\smallskip
\noindent {\rm b)} $w=\cdot 1$: We already know from Proposition \ref{prop:D-numbers} that $R_{\cdot 1}= V(3,1,1)$.
\smallskip

\noindent {\rm c)} $w=\cdot 00$: Using  the Propagation Principle, we see that a digit block $\cdot 10$ is always followed directly by the first element of a trident of $\cdot 00$'s and vice versa. This implies the statement in {\rm c)}, because of {\rm b)}.

\smallskip

\noindent {\rm d)} $w=\cdot 01$: This result is given in Remark 6.2 in the paper \cite{Dekk-phi-FQ}.

\smallskip

\noindent {\rm e)} $w=\cdot 000$: Using  the Propagation Principle, we see that a $\cdot 100$ is always followed directly by the first element of a trident of $\cdot 000$'s and vice versa. So {\rm e)} is implied by {\rm g)}.

\smallskip

\noindent {\rm f)} $w=\cdot 001$: Take the first sequence $V(3,1,2)$ of $R_{\cdot 00}$, and put $p=3, q=1, r=2$. Then the first sequence of   $R_{\cdot 000}$ is equal to $V(4,3,-1)=V(p+q,p,r-p)$. It then follows from Lemma \ref{lem:VAconv} that the first sequence of $R_{\cdot 001}$ is equal to $V(2p+q,p+q,r)=V(7,4,2)$.

\smallskip

\noindent {\rm g)} $w=\cdot 100$:   For the first 17 numbers we check that $\cdot 100$ occurs as prefix of $\beta^-(N)$ if and only if $1000$ occurs as suffix of  $\beta^+(N)$.
   The result then follows from Theorem \ref{th:phi-w0}: $R_w= V(L_{m-1},L_{m-2},\gamma_w) \;\;\text{\rm if}\; w_{m-1}=0$, where here $m=4$, so
$R_{1000} = V(L_{3},L_{2},\gamma_{1000})=V(4,3,-2)$. Here $\gamma_{1000}$  is determined by noting that $N=5$ is the first number in $R_{1000}$.

\smallskip 

\noindent {\rm h)} $w= \cdot 101$: Take the sequence $R_{\cdot 10}=V(3,1,1)$, and put $p=3, q=1, r=1$. Then $R_{\cdot 100}$ is equal to $V(4,3,-2)=V(p+q,p,r-p)$. It then follows from Lemma \ref{lem:VAconv} that the sequence  $R_{\cdot 101}$ is equal to $V(2p+q,p+q,r)=V(7,4,1)$.
\hfill $\Box$

\medskip

The reader might think that we can now proceed, as we did earlier, from these cases to words $w$ with larger lengths $m$, using the same tools. However, this does not work. The reason is that the $\beta^-(N)$ words do not occur in lexicographical order, in contrast with the $\beta^+(N)$ words.
Some occurrence sequences are Lucas-Wythoff, some are not---but still close to Lucas-Wythoff sequences.

Recall the  three (Sturmian) morphisms $f,g$ and $h$ from Equation (\ref{eq:Fib3}).
Note that $f$ equals the square of the Fibonacci morphism $a\mapsto ab, \, b\mapsto a$, so $f$ has fixed point $\xF$, the Fibonacci word.
The fixed points $\xG, \xH$ of $g$ and $h$ are given by $\xG=b\,\xF,\,\xH=a\,\xF$ ---see \cite{Berstel-Seebold} Theorem 3.1.

Let $V_{\F}, V_{\G} ,V_{\Ha}$ denote the families of sequences  having $\xF, \xG, \xH$ as first differences, with first element an arbitrary integer. Then, by definition, one example  is $V=V_{\F}$, if we take $V_{\F}(1)=p+q+r$. We also already have encountered  an $V_{\G}$, since $V_0=V_{\G}$, if we take  $V_{\G}(1)=r$. This follows from $V_0(p,q,r)= r, p+q+r, \dots = r, b+r, \dots$, which gives $\Delta V_0 = b \xF =\xG$. We mention that one can show that there do not exist  $\alpha, p, q$, and $r$ such that $V_{\Ha}$ is a generalized Beatty sequence  $V = (p\lfloor n \alpha \rfloor + q n +r) $.

We conjecture that the following holds.

\begin{conjecture*} Let $\beta(N)=\beta^+(N)\cdot\beta^-(N)$ be the base phi expansion of the   number $N$.\\
  Let $w$ be a word of length $m$. Let $R_{\cdot w}$ be the sequence of occurrences  of numbers $N$ such that  the  first $m$ digits of $\beta^-(N)$  are equal to $w$, i.e., $d_{-1}\dots d_{-m}(N)=w$. Then  there exist two Lucas numbers $a$ and $b$ such that either $R_{\cdot w} = V_{\F},$ or $ R_{\cdot w} = V_{\G},$ or $R_{\cdot w} = V_{\Ha}$. A second possibility is that $R_{\cdot w}$ is a union of three of such sequences.
\end{conjecture*}

\medskip

In all  cases  in Proposition \ref{prop:smallminus} the sequence  $ R_{\cdot w}$ is a  $V_{\F}$, except $ R_{\cdot 010}$, which is a union of three $V_{\G}$'s, the middle one being $V_{\G}(4,3,-4)$. The first case where a $V_{\G}$ as $ R_{\cdot w}$ occurs, is for $w=\cdot 1001$, where $a=29, b=18$. The first case where $V_{\Ha}$ as a $R_{\cdot w}$   occurs, is as first element of the trident for the digit block $w=\cdot 0100$, where $a=18, b=11$.

\noindent AMS Classification Numbers: 11D85, 11A63, 11B39

\end{document}